\documentclass[11pt]{article}
\usepackage{amsfonts}

\usepackage[dvips]{graphicx}
\usepackage{amsmath}

\newtheorem{theorem}{Theorem}[subsection]

\newtheorem{corollary}[theorem]{Corollary}

\newtheorem{example}[theorem]{Example}

\newtheorem{lemma}[theorem]{Lemma}

\newtheorem{proposition}[theorem]{Proposition}

\newenvironment{proof}[1][Proof]{\textbf{#1.} }{\ \rule{0.5em}{0.5em}}

\input{tcilatex}

\begin{document}

\title{Universal graded characters and limit of Lusztig $q$-analogues}
\author{C\'{e}dric Lecouvey \\
Laboratoire de Math\'{e}matiques Pures et Appliqu\'{e}es Joseph Liouville\\
B.P. 699 62228 Calais Cedex\\
Cedric.Lecouvey@lmpa.univ-littoral.fr}
\date{}
\maketitle

\begin{abstract}
Let $G$ be a symplectic or orthogonal complex Lie group with Lie algebra $%
\frak{g}$. As a $G$-module, the decomposition of the symmetric algebra $S(%
\frak{g)}$ into its irreducible components can be explicitely obtained by
using identities due to Littlewood. We show that the multiplicities
appearing in the decomposition of the $k$-th graded component of $S(\frak{g})
$ do not depend on the rank $n$ of $\frak{g}$ \ providing $n$ is
sufficiently large. Thanks to a classical result by Kostant, we establish a
similar result for the $k$-th graded component of the space $H(\frak{g})$ of 
$G$-harmonic polynomials. These stabilization properties are equivalent to
the existence of a limit in infinitely many variables for the graded
characters associated to $S(\frak{g})$ and $H(\frak{g})$. The limits so
obtained are formal series with coefficients in the ring of universal
characters introduced by Koike and Terada.

From Hesselink expression of the graded character of harmonics, the
coefficient of degree $k$ in the Lusztig $q$-analogue $K_{\lambda ,\emptyset
}^{\frak{g}}(q)$ associated to the fixed partition $\lambda $ thus stablizes
for $n$ sufficiently large. By using Morris-type recurrence formulas, we
prove that this is also true for the polynomials $K_{\lambda ,\mu }^{\frak{g}%
}(q)$ where $\mu $ is a nonempty fixed partition. This can be reformulated
in terms of a stability property for the dimension of the components of the
Brylinski-Kostant filtration. We also associate to each pair of partitions $%
(\lambda ,\mu )$ formal series $K_{\lambda ,\mu }^{\frak{so}}(q)$ and $%
K_{\lambda ,\mu }^{\frak{sp}}(q)$, which can be regarded as natural limit of
the Lusztig $q$-analogues.\ One gives a duality property for these limits
and obtains simple expressions when $\lambda $ is a row or a column
partition.
\end{abstract}

\section{Introduction}

The multiplicity $K_{\lambda ,\mu }$ of the weight $\mu $ in the irreducible
finite-dimensional representation $V^{\frak{g}}(\lambda )$ of the simple Lie
group $G$ with Lie algebra $\frak{g}$ can be written in terms of the
ordinary Kostant partition function $\mathcal{P}$ defined by the equality 
\begin{equation*}
\prod_{\alpha \text{ positive root}}\dfrac{1}{(1-e^{\alpha })}=\sum_{\beta }%
\mathcal{P}(\beta )e^{\beta }
\end{equation*}
where $\beta $ runs on the set of nonnegative integral combinations of
positive roots of $\frak{g}$. Thus $\mathcal{P}(\beta )$ is the number of
ways the weight $\beta $ can be expressed as a sum of positive roots. Then,
one derives from the Weyl character formula 
\begin{equation}
K_{\lambda ,\mu }^{\frak{g}}=\sum_{w\in W^{\frak{g}}}(-1)^{\ell (w)}\mathcal{%
P}(w(\lambda +\rho )-(\mu +\rho ))  \label{Kosn}
\end{equation}
where $W^{\frak{g}}$ is the Weyl group of $\frak{g}$.

\noindent The Lusztig $q$-analogue of weight multiplicity $K_{\lambda ,\mu
}^{\frak{g}}(q)$ is obtained by substituting the ordinary Kostant partition
function $\mathcal{P}$ by its $q$-analogue $\mathcal{P}_{q}$ in (\ref{Kosn}%
).\ Namely $\mathcal{P}_{q}$ is defined by the equality 
\begin{equation*}
\prod_{\alpha \text{ positive root}}\dfrac{1}{(1-qe^{\alpha })}=\sum_{\beta }%
\mathcal{P}_{q}(\beta )e^{\beta }
\end{equation*}
and we have 
\begin{equation*}
K_{\lambda ,\mu }^{\frak{g}}(q)=\sum_{w\in W^{\frak{g}}}(-1)^{\ell (w)}%
\mathcal{P}_{q}(w(\lambda +\rho )-(\mu +\rho )).
\end{equation*}
As shown by Lusztig \cite{Lut}, $K_{\lambda ,\mu }^{\frak{g}}(q)$ is a
polynomial in $q$ with nonnegative integer coefficients. Many
interpretations of the Lusztig $q$-analogues exist. In particular, they can
be obtained from the Brylinski-Kostant filtration of weight spaces \cite{bry}%
.\ The polynomials $K_{\lambda ,\emptyset }^{\frak{g}}(q)$ appear in the
graded character of the harmonic polynomials associated to $\frak{g}$ \cite
{hes}. We also recover the Lusztig $q$-analogues as the coefficients of the
expansion of the Hall-Littlewood polynomials on the basis of Weyl characters
(see \cite{NR}). This notably permits to prove that there are affine
Kazhdan-Lusztig polynomials. In \cite{LSc1} Lascoux and Sch\"{u}tzenberger
have obtained a combinatorial expression for $K_{\lambda ,\mu }^{\frak{gl}%
_{n}}(q)$ in terms of the charge statistic on the semistandard tableaux of
shape $\lambda $ and evaluation $\mu .$ By using the combinatorics of
crystal graphs introduced by Kashiwara and Nakashima \cite{KN}, we have also
established similar formulas \cite{lec}, \cite{lec2} for the Lusztig $q$%
-analogues associated to the symplectic and orthogonal Lie algebras when $%
(\lambda ,\mu )$ satisfies restrictive constraints.

\bigskip

Consider $\lambda ,\mu $ two partitions of length at most $m.\;$These
partitions can be regarded as dominant weights for $\frak{g=gl}_{n},\frak{so}%
_{2n+1},\frak{sp}_{2n}$ or $\frak{so}_{2n}$ when $n\geq m.$ Then, $%
K_{\lambda ,\mu }^{\frak{gl}_{n}}(q)$ does not depend on the rank $n$
considered. Such a property does not hold for the Lusztig $q$-analogues $%
K_{\lambda ,\mu }^{\frak{g}}(q)$ when $\frak{g=so}_{2n+1},\frak{sp}_{2n}$ or 
$\frak{so}_{2n}$ which depend in general on the rank of the Lie algebra
considered. Write 
\begin{equation*}
K_{\lambda ,\mu }^{\frak{g}}(q)=\sum_{k\geq 0}K_{\lambda ,\mu }^{\frak{g}%
,k}q^{k}.
\end{equation*}
We first establish in this paper that for $\frak{g=so}_{2n+1},\frak{sp}_{2n}$
or $\frak{so}_{2n},$ the coefficient $K_{\lambda ,\mu }^{\frak{g},k}$
stabilizes when $n$ tends to the infinity.\ More precisely, $K_{\lambda ,\mu
}^{\frak{g},k}$ does not depend on the rank $n$ of $\frak{g}$ providing $%
n\geq 2k+a$ where $a$ is the number of nonzero parts of $\mu $ (Theorem \ref
{Th-Stab}).\ By Brylinski's interpretation of the coefficients $K_{\lambda
,\mu }^{\frak{g},k}$ \cite{bry}, one then obtains that the dimension of the $%
k$-th component of the Brylinski-Kostant filtration associated to the
finite-dimensional irreducible representations of $\frak{g=so}_{2n+1},\frak{%
sp}_{2n}$ or $\frak{so}_{2n}$ stabilizes for $n$ sufficiently large (Theorem 
\ref{th_stabbry}). Observe that this stabilization is immediate for $\frak{%
g=gl}_{n}$ since the polynomials $K_{\lambda ,\mu }^{\frak{gl}_{n}}(q)$ does
not depend on $n$. For $\frak{g=so}_{2n+1},\frak{sp}_{2n}$ or $\frak{so}_{2n}
$ the Brylinski-Kostant filtration depends in general on the rank considered
and it seems difficult to obtain the dimension of its components by direct
computations.

\noindent Our methods is as follows.\ We obtain the explicit decomposition
of the symmetric algebra $S(\frak{g)}$ considered as a $G$-module into its
irreducible components by using identities due to Littlewood. This permits
to show that the multiplicities appearing in the decomposition of the $k$-th
graded component of $S(\frak{g})$ do not depend on the rank $n$ of $\frak{g}$
providing $n$ is sufficiently large. Thanks to a classical result by
Kostant, we establish a similar result for the $k$-th graded component of
the space $H(\frak{g})$ of $G$-harmonic polynomials. These stabilization
properties is equivalent to the existence of a limit in infinitely many
variables for the graded characters associated to $S(\frak{g})$ and $H(\frak{%
g})$. The limits so obtained are formal series with coefficients in the ring
of universal characters introduced by Koike and Terada. From Hesselink
expression \cite{hes} of the graded character of $H(\frak{g})$, one then
derives that $K_{\lambda ,\emptyset }^{\frak{g},k}$ stabilizes for $n$
sufficiently large. By using Morris-type recurrence formulas for the Lusztig 
$q$-analogues \cite{lec2}, we prove that this is also true for the
coefficients $K_{\lambda ,\mu }^{\frak{g},k}$ where $\mu $ is a nonempty
fixed partition. We also observe that these formulas permit to give an
explicit lower bound for the degree of the $q$-analogues $K_{\lambda ,\mu }^{%
\frak{g},}(q)$ such that $K_{\lambda ,\mu }^{\frak{g},}(q)\neq 0$.\ We
establish that the limits of the coefficients $K_{\lambda ,\mu }^{\frak{so}%
_{2n+1},k}$ and $K_{\lambda ,\mu }^{\frak{so}_{2n+1},k}$ are the same. Write 
$K_{\lambda ,\mu }^{\frak{sp},k}$ and $K_{\lambda ,\mu }^{\frak{so},k}$
respectively for the limits of the coefficients $K_{\lambda ,\mu }^{\frak{so}%
_{2n+1},k}$ and $K_{\lambda ,\mu }^{\frak{sp}_{2n},k}$ when $n$ tends to the
infinity.

\noindent The stabilization property of the coefficients $K_{\lambda ,\mu }^{%
\frak{g},k}$ suggests then to introduce the formal series 
\begin{equation*}
K_{\lambda ,\mu }^{\frak{so}}(q)=\sum_{k\geq 0}K_{\lambda ,\mu }^{\frak{so}%
,k}q^{k}\text{ and }K_{\lambda ,\mu }^{\frak{sp}}(q)=\sum_{k\geq
0}K_{\lambda ,\mu }^{\frak{sp},k}q^{k}.
\end{equation*}
These series belong $\mathbb{N}[[q]]$ and can be regarded as natural limits
of the polynomials $K_{\lambda ,\mu }^{\frak{g}}(q).$ We establish a duality
result between the formal series $K_{\lambda ,\emptyset }^{\frak{so}}(q)$
and $K_{\lambda ,\emptyset }^{\frak{sp}}(q)$ (Theorem \ref{th_dual}).
Namely, we have 
\begin{equation}
K_{\lambda ,\emptyset }^{\frak{so}}(q)=K_{\lambda ^{\prime },\emptyset }^{%
\frak{sp}}(q)  \label{conjBC}
\end{equation}
where $\lambda ^{\prime }$ is the conjugate partition of $\lambda $. Note
that (\ref{conjBC}) do not hold in general if we replace the formal series $%
K_{\lambda ,\mu }^{\frak{so}}(q)$ and $K_{\lambda ,\mu }^{\frak{sop}}(q)$ by
the polynomials $K_{\lambda ,\mu }^{\frak{g},k}(q).$ We also give recurrence
formulas (\ref{rec_for_limit}), (\ref{rec_for_limit_0}) for the series $%
K_{\lambda ,\mu }^{\frak{so}}(q)$ and $K_{\lambda ,\mu }^{\frak{sop}}(q)$.\
Thanks to these recurrence formulas, one derives simple expressions for the
formal series $K_{\lambda ,\emptyset }^{X}(q)$ when $\lambda $ is a row or a
column partition (Proposition \ref{prop_expli}). Note that we have not find
so simple formulas for the Lusztig $q$-analogues $K_{\lambda ,\mu }^{\frak{g}%
}(q)$ even in the cases when $\lambda $ is a column or a row partition.
Moreover the duality (\ref{conjBC}) is false in general for the polynomials $%
K_{\lambda ,\mu }^{\frak{g}}(q)$. This suggests that the study of the series 
$K_{\lambda ,\mu }^{\frak{so}}(q)$ and $K_{\lambda ,\mu }^{\frak{sp}}(q)$
which is initiated in this paper, could be easier than that of the Lusztig $q
$-analogues.

\bigskip

The paper is organized as follows.\ In Section $2$ we recall the necessary
background on symplectic and orthogonal Lie algebras, universal characters,
and Lusztig $q$-analogues which is needed in the sequel. In Section $3$ we
introduce universal graded characters as limit in infinitely many variables
for the graded characters associated to $S(\frak{g})$ and $H(\frak{g})$. We
obtain the stabilization property of the coefficients $K_{\lambda ,\mu }^{%
\frak{g},k}$ in Section $4$ and reformulate this result in terms of the
Brylinski-Kostant filtration. In Section $5,$ we introduce the formal series 
$K_{\lambda ,\mu }^{\frak{so}}(q)$ and $K_{\lambda ,\mu }^{\frak{sop}}(q)$,
establish recurrence formulas which permit to compute them by induction,
prove the duality (\ref{conjBC}) and give explicit formulas for $K_{\lambda
,\mu }^{\frak{so}}(q)$ and $K_{\lambda ,\mu }^{\frak{sop}}(q)$ when $\lambda 
$ is a row or a column partition. Finally in Section $6$ we have added a few
considerations on the possibility to define Hall-Littlewood polynomials from
the formal series $K_{\lambda ,\mu }^{\frak{so}}(q)$ and $K_{\lambda ,\mu }^{%
\frak{sop}}(q)$.

\section{Background}

\subsection{Convention for the root systems of types $B,C$ and $D$}

In the sequel $G$ is one of the complex Lie groups $Sp_{2n},SO_{2n+1}$ or $%
SO_{2n}$ and $\frak{g}$ its Lie algebra.\ We follow the convention of \cite
{KT} to realize $G$ as a subgroup of $GL_{N}$ and $\frak{g}$\ as a
subalgebra of $\frak{gl}_{N}$ where 
\begin{equation*}
N=\left\{ 
\begin{tabular}{l}
$2n$ when $G=Sp_{2n}$ \\ 
$2n+1$ when $G=SO_{2n+1}$ \\ 
$2n$ when $G=SO_{2n}$%
\end{tabular}
\right. .
\end{equation*}
With this convention the maximal torus $T$ of $G$ and the Cartan subalgebra $%
\frak{h}$ of $\frak{g}$ coincide respectively with the subgroup and the
subalgebra of diagonal matrices of $G$ and $\frak{g}$. Similarly the Borel
subgroup $B$ of $G$ and the Borel subalgebra $\frak{b}_{+}$ of $\frak{g}$
coincide respectively with the subgroup and subalgebra of upper triangular
matrices of $G$ and $\frak{g}$. This gives the triangular decomposition $%
\frak{g=b}_{+}\oplus \frak{h}\oplus \frak{b}_{-}$ for the Lie algebra $\frak{%
g}$. Let $e_{i},h_{i},f_{i},$ $i\in \{1,...,n\}$ be a set of Chevalley
generators such that $e_{i}\in \frak{b}_{+},h_{i}\in \frak{h}$ and $f_{i}\in 
\frak{b}_{-}$ for any $i$.

\noindent Let $d_{N}$ be the linear subspace of $\frak{gl}_{N}$ consisting
of the diagonal matrices.\ For any $i\in\{1,...,n\},$ write $\varepsilon_{i}$
for the linear map $\varepsilon_{i}:d_{N}\rightarrow\mathbb{C}$ such that $%
\varepsilon_{i}(D)=\delta_{i}$ for any diagonal matrix $D$ whose $(i,i)$%
-coefficient is $\delta_{i}.$ Then $(\varepsilon_{1},...,\varepsilon _{n})$
is an orthonormal basis of the Euclidean space $\frak{h}_{\mathbb{R}}^{\ast}$
(the real part of $\frak{h}^{\ast}).$ We denote by $<\cdot,\cdot>$ the usual
scalar product on $\frak{h}_{\mathbb{R}}^{\ast}$.

\bigskip

\noindent Let $R$ be the root system associated to $G.\;$We can take for the
simple roots of $\frak{g}$%
\begin{equation}
\left\{ 
\begin{tabular}{l}
$\Sigma ^{+}=\{\alpha _{n}=\varepsilon _{n}\text{ and }\alpha
_{i}=\varepsilon _{i}-\varepsilon _{i+1}\text{, }i=1,...,n-1\text{ for the
root system }B_{n}\}$ \\ 
$\Sigma ^{+}=\{\alpha _{n}=2\varepsilon _{n}\text{ and }\alpha
_{i}=\varepsilon _{i}-\varepsilon _{i+1}\text{, }i=1,...,n-1\text{ for the
root system }C_{n}\}$ \\ 
$\Sigma ^{+}=\{\alpha _{n}=\varepsilon _{n}+\varepsilon _{n-1}\text{ and }%
\alpha _{i}=\varepsilon _{i}-\varepsilon _{i+1}\text{, }i=1,...,n-1\text{
for the root system }D_{n}\}$%
\end{tabular}
\right. .  \label{simple_roots}
\end{equation}
Then the sets of positive roots are 
\begin{equation*}
\left\{ 
\begin{tabular}{l}
$R^{+}=\{\varepsilon _{i}-\varepsilon _{j},\varepsilon _{i}+\varepsilon _{j}%
\text{ with }1\leq i<j\leq n\}\cup \{\varepsilon _{i}\text{ with }1\leq
i\leq n\}\text{ for the root system }B_{n}$ \\ 
$R^{+}=\{\varepsilon _{i}-\varepsilon _{j},\varepsilon _{i}+\varepsilon _{j}%
\text{ with }1\leq i<j\leq n\}\cup \{2\varepsilon _{i}\text{ with }1\leq
i\leq n\}\text{ for the root system }C_{n}$ \\ 
$R^{+}=\{\varepsilon _{i}-\varepsilon _{j},\varepsilon _{i}+\varepsilon _{j}%
\text{ with }1\leq i<j\leq n\}\text{ for the root system }D_{n}$%
\end{tabular}
\right. .
\end{equation*}
We denote by $R$ the set of roots of $G.\;$For any $\alpha \in R$, let $%
\alpha ^{\vee }=\frac{\alpha }{<\alpha ,\alpha >}$ be the coroot
corresponding to $\alpha $. The Weyl group of the Lie group $G$ is the
subgroup of the permutation group of the set $\{\overline{n},...,\overline{2}%
,\overline{1},1,2,...,n\}$\ generated by the permutations 
\begin{equation*}
\left\{ 
\begin{tabular}{l}
$s_{i}=(i,i+1)(\overline{i},\overline{i+1}),$ $i=1,...,n-1$ and $s_{n}=(n,%
\overline{n})$ $\text{for the root systems }B_{n}$ and $C_{n}$ \\ 
$s_{i}=(i,i+1)(\overline{i},\overline{i+1}),$ $i=1,...,n-1$ and $%
s_{n}^{\prime }=(n,\overline{n-1})(n-1,\overline{n})$ $\text{for the root
system }D_{n}$%
\end{tabular}
\right.
\end{equation*}
where for $a\neq b$ $(a,b)$ is the simple transposition which switches $a$
and $b.$ We identify the subgroup of $W^{\frak{g}}$ generated by $%
s_{i}=(i,i+1)(\overline{i},\overline{i+1}),$ $i=1,...,n-1$ with the
symmetric group $S_{n}.$ We denote by $\ell $ the length function
corresponding to the above set of generators. The action of $w\in W^{\frak{g}%
}$ on $\beta =(\beta _{1},...,\beta _{n})\in \frak{h}_{\mathbb{R}}^{\ast }$
is defined by 
\begin{equation*}
w\cdot (\beta _{1},...,\beta _{n})=(\beta _{1}^{w^{-1}},...,\beta
_{n}^{w^{-1}})
\end{equation*}
where $\beta _{i}^{w}=\beta _{w(i)}$ if $w(i)\in \{1,...,n\}$ and $\beta
_{i}^{w}=-\beta _{w(\overline{i})}$ otherwise. We denote by $\rho $ the half
sum of the positive roots of $R^{+}$.\ The dot action of $W^{\frak{g}}$ on $%
\beta =(\beta _{1},...,\beta _{n})\in \frak{h}_{\mathbb{R}}^{\ast }$ is
defined by 
\begin{equation}
w\circ \beta =w\cdot (\beta +\rho )-\rho .  \label{dotaction}
\end{equation}
Write $P$ and $P^{+}$ for the weight lattice and the cone of dominant
weights of $G.\;$As usual we consider the order on $P$ defined by $\beta
\leq \gamma $ if and only if $\beta -\gamma \in Q^{+}.$

\noindent For any positive integer $m$, denote by $\mathcal{P}_{m}$ the set
of partitions with at most $m$ nonzero parts. Let $\mathcal{P}_{m}(k),k\in 
\mathbb{N}$ be the subset of $\mathcal{P}_{m}$ containing the partitions $%
\lambda $ such that $\left| \lambda \right| =\lambda _{1}+\cdot \cdot \cdot
+\lambda _{m}=k.$ Set $\mathcal{P=\cup }_{m\in \mathbb{N}}\mathcal{P}_{m}$
and $\mathcal{P}_{m}[k]=\mathcal{\cup }_{a\leq k}\mathcal{P}_{m}(a)$.

\noindent Each partition $\lambda =(\lambda _{1},...,\lambda _{n})\in $ $%
\mathcal{P}_{n}$ will be identified with the dominant weight $%
\sum_{i=1}^{n}\lambda _{i}\varepsilon _{i}.\;$Then the irreducible
finite-dimensional polynomial representations of $G$ are parametrized by the
partitions of $\mathcal{P}_{n}$.\ For any $\lambda \in \mathcal{P}_{n},$
denote by $V^{\frak{g}}(\lambda )$ the irreducible finite-dimensional
representation of $G$ of highest weight $\lambda .$ The representation $V^{%
\frak{g}}(1)$ associated to the partition $\lambda =(1)$ is called the
vector representation of $G$. For any weight $\beta \in P$ and any partition 
$\lambda \in \mathcal{P}_{n},$ we write $V^{\frak{g}}(\lambda )_{\beta }$
for the weight space associated to $\beta $ in $V^{\frak{g}}(\lambda )$.

\noindent We denote by $Q$ the root lattice of $\frak{g}$ and write $Q^{+}$
for the elements of $Q$ which are linear combination of positive roots with
nonnegative coefficients.

\noindent The exponents $\{m_{1},...,m_{n}\}$ of the root system $R$
verifies $m_{i}=2i-1,$ $i=1,...,n$ when $R$ is of type $B_{n}$ or $C_{n}$
and 
\begin{equation}
m_{i}=2i-1,i\in\{1,...,n-1\}\text{ and }m_{n}=n-1  \label{def_expo}
\end{equation}
when $R$ is of type $D_{n}.\;$

\bigskip

\noindent \textbf{Remarks:}

\noindent$\mathrm{(i):}$ The integer $n-1$ appears twice in the exponents of
a root system of type $D_{n}$ when $n$ is even.

\noindent$\mathrm{(ii):}$ The exponents $m_{i},i=1,...,n-1$ are the same for
the three root systems of type $B_{n},C_{n}$ or $D_{n}.$

\bigskip

\noindent As customary, we identify $P$ the lattice of weights of $G$ with a
sublattice of $(\frac{1}{2}\mathbb{Z})^{n}.$ For any $\beta=(\beta
_{1},...,\beta_{n})\in P,$ we set $\left| \beta\right| =\beta_{1}+\cdot
\cdot\cdot+\beta_{n}.\;$We use for a basis of the group algebra $\mathbb{Z}[%
\mathbb{Z}^{n}],$ the formal exponentials $(e^{\beta})_{\beta\in \mathbb{Z}%
^{n}}$ satisfying the relations $e^{\beta_{1}}e^{\beta_{2}}=e^{\beta_{1}+%
\beta_{2}}.$ We furthermore introduce $n$ independent indeterminates $%
x_{1},...,x_{n}$ in order to identify $\mathbb{Z}[\mathbb{Z}^{n}]$ with the
ring of polynomials $\mathbb{Z}[x_{1},...,x_{n},x_{1}^{-1},...,x_{n}^{-1}]$
by writing $e^{\beta}=x_{1}^{\beta_{1}}\cdot\cdot\cdot
x_{n}^{\beta_{n}}=x^{\beta}$ for any $\beta=(\beta _{1},...,\beta_{n})\in%
\mathbb{Z}^{n}.$

\noindent Write $s_{\lambda }^{\frak{gl}_{n}}$ for the Weyl character (Schur
function) of $V^{\frak{gl}_{n}}(\lambda )$ the finite-dimensional $\frak{gl}%
_{n}$-module of highest weight $\lambda .$ The character ring of $GL_{n}$ is 
$\Lambda _{n}=\mathbb{Z[}x_{1},...,x_{n}]^{\mathrm{sym}}$ the ring of
symmetric functions in $n$ variables.\ 

\noindent For any $\lambda \in \mathcal{P}_{n},$ we denote by $s_{\lambda }^{%
\frak{g}}$ the Weyl character of $V^{\frak{g}}(\lambda )$. Let $R^{\frak{g}}$
be the character ring of $G.$ Then 
\begin{equation*}
R^{\frak{g}}=\mathbb{Z}[x_{1},...,x_{n},x_{1}^{-1},...,x_{n}^{-1}]^{W^{\frak{%
g}}}
\end{equation*}
is the $\mathbb{Z}$-algebra with basis $\{s_{\lambda }^{\frak{g}}\mid
\lambda \in \mathcal{P}_{n}\}.$

\noindent In the sequel we will suppose $n\geq 2$ when $\frak{g}=\frak{sp}%
_{2n}$ or $\frak{so}_{2n+1}$ and $n\geq 4$ when $\frak{g}=\frak{so}_{2n}.$

\bigskip

\noindent For each Lie algebra $\frak{g=so}_{N}$ or $\frak{sp}_{N}$ and any
partition $\nu \in \mathcal{P}_{N}$, we denote by $V^{\frak{gl}_{N}}(\nu
)\downarrow _{\frak{g}}^{\frak{gl}_{N}}$ the restriction of $V^{\frak{gl}%
_{N}}(\nu )$ to $\frak{g}$. Set 
\begin{eqnarray}
V^{\frak{gl}_{N}}(\nu ) &\downarrow &_{\frak{g}}^{\frak{gl}%
_{N}}=\bigoplus_{\lambda \in \mathcal{P}_{n}}V^{\frak{so}_{N}}(\lambda
)^{\oplus b_{\nu ,\lambda }^{\frak{so}_{N}}},  \label{dec_char} \\
V^{\frak{gl}_{2n}}(\nu ) &\downarrow &_{\frak{g}}^{\frak{gl}%
_{2n}}=\bigoplus_{\lambda \in \mathcal{P}_{n}}V^{\frak{so}_{2n}}(\lambda
)^{\oplus b_{\nu ,\lambda }^{\frak{sp}_{2n}}}.  \notag
\end{eqnarray}
This define in particular the branching coefficients $b_{\nu ,\lambda }^{%
\frak{so}_{N}}$ and $b_{\nu ,\lambda }^{\frak{sp}_{2n}}$.\ The restriction
map $r^{\frak{g}}$ is defined by setting 
\begin{equation*}
r^{\frak{g}}:\left\{ 
\begin{array}{c}
\mathbb{Z[}x_{1},...,x_{N}]^{\mathrm{sym}}\rightarrow R^{\frak{g}} \\ 
s_{\nu }^{\frak{gl}_{N}}\longmapsto \mathrm{char}(V^{\frak{gl}_{N}}(\nu
)\downarrow _{\frak{g}}^{\frak{gl}_{N}})
\end{array}
\right. .
\end{equation*}
We have then 
\begin{equation*}
r^{\frak{g}}(s_{\nu }^{\frak{gl}_{N}})=\left\{ 
\begin{array}{c}
s_{\nu }^{\frak{gl}_{N}}(x_{1},...,x_{n},x_{n}^{-1},...,x_{1}^{-1})\text{
when }N=2n \\ 
s_{\nu }^{\frak{gl}_{N}}(x_{1},...,x_{n},0,x_{n}^{-1},...,x_{1}^{-1})\text{
when }N=2n+1
\end{array}
\right. .
\end{equation*}
Let $\mathcal{P}_{n}^{(2)}$ and $\mathcal{P}_{n}^{(1,1)}$ be the subsets of $%
\mathcal{P}_{n}$ containing the partitions with even length rows and the
partitions with even length columns, respectively. When $\nu \in \mathcal{P}%
_{n}$ we have the following formulas for the branching coefficients $b_{\nu
,\lambda }^{\frak{so}_{N}}$ and $b_{\nu ,\lambda }^{\frak{p}_{2n}}$ :

\begin{proposition}
\label{prop_lit}(see \cite{Li} appendix p\ 295)\label{prop_in_A2n}

\noindent Consider $\nu \in \mathcal{P}_{n}.$ Then:

\begin{enumerate}
\item  $b_{\nu ,\lambda }^{\frak{so}_{2n+1}}=b_{\nu ,\lambda }^{\frak{so}%
_{2n}}=\sum_{\gamma \in \mathcal{P}_{n}^{(2)}}c_{\lambda ,\gamma }^{\nu }$

\item  $b_{\nu ,\lambda }^{\frak{sp}_{2n}}=\sum_{\gamma \in \mathcal{P}%
_{n}^{(1,1)}}c_{\lambda ,\gamma }^{\nu }$
\end{enumerate}

\noindent where $c_{\gamma ,\lambda }^{\nu }$ is the usual
Littlewood-Richardson coefficient corresponding to the partitions $\gamma
,\lambda $ and $\nu $.
\end{proposition}

\noindent Note that the equality $b_{\nu ,\lambda }^{\frak{so}%
_{2n+1}}=b_{\nu ,\lambda }^{\frak{so}_{2n}}$ becomes false in general when $%
\nu \notin \mathcal{P}_{n}$.

\noindent As suggested by Proposition \ref{prop_lit}, the manipulation of
the Weyl characters is simplified by working with infinitely many
variables.\ In \cite{K}, Koike and Terada have introduced a universal
character ring for the classical Lie groups.\ This ring can be regarded as
the ring $\Lambda =\mathbb{Z[}x_{1},...,x_{n},...]^{\mathrm{sym}}$ of
symmetric functions in countably many variables.\ It is equipped with three
natural $\mathbb{Z}$-bases indexed by partitions, namely 
\begin{equation}
\mathcal{B}^{\frak{gl}}\mathcal{=}\{\mathtt{s}_{\lambda }^{\frak{gl}}\mid
\lambda \in \mathcal{P}\},\text{ }\mathcal{B}^{\frak{sp}}\mathcal{=}\{%
\mathtt{s}_{\lambda }^{\frak{sp}}\mid \lambda \in \mathcal{P}\}\text{ and }%
\mathcal{B}^{\frak{so}}\mathcal{=}\{\mathtt{s}_{\lambda }^{\frak{so}}\mid
\lambda \in \mathcal{P}\}.  \label{bases}
\end{equation}
We have in particular the decompositions: 
\begin{equation}
\mathtt{s}_{\nu }^{\frak{gl}}=\sum_{\lambda \in \mathcal{P}}\sum_{\gamma \in 
\mathcal{P}_{n}^{(2)}}c_{\lambda ,\gamma }^{\nu }\mathtt{s}_{\lambda }^{%
\frak{so}}\text{ and }\mathtt{s}_{\nu }^{\frak{gl}}=\sum_{\lambda \in 
\mathcal{P}}\sum_{\gamma \in \mathcal{P}_{n}^{(1,1)}}c_{\lambda ,\gamma
}^{\nu }\mathtt{s}_{\lambda }^{\frak{sp}}.  \label{deuSgl}
\end{equation}
We denote by $\varphi $ the linear involution defined on $\Lambda $ by 
\begin{equation}
\varphi (\mathtt{s}_{\lambda }^{\frak{so}})=\mathtt{s}_{\lambda ^{\prime }}^{%
\frak{sp}}.  \label{def_fi}
\end{equation}
For any positive integer $n,$ denote by $\Lambda _{n}=\mathbb{Z[}%
x_{1},...,x_{n}]^{\mathrm{sym}}$ the ring of symmetric functions in $n$
variables.\ Write 
\begin{equation*}
\pi _{n}:\mathbb{Z[}x_{1},...,x_{n},...]^{\mathrm{sym}}\rightarrow \mathbb{Z[%
}x_{1},...,x_{n}]^{\mathrm{sym}}
\end{equation*}
for the ring homomorphism obtained by specializing each variable $x_{i},i>n$
at $0.$ Then $\pi _{n}(\mathtt{s}_{\lambda }^{\frak{gl}})=s_{\lambda }^{%
\frak{gl}_{n}}.$ Let $\pi ^{\frak{sp}_{2n}}$ and $\pi ^{\frak{so}_{N}}$ be
the specialization homomorphisms defined by setting $\pi ^{\frak{sp}%
_{2n}}=r^{\frak{sp}_{2n}}\circ \pi _{2n}$ and $\pi ^{\frak{so}_{N}}=r^{\frak{%
so}_{N}}\circ \pi _{N}$.\ For any partition $\lambda \in \mathcal{P}_{n}$
one has 
\begin{equation*}
s_{\lambda }^{\frak{sp}_{2n}}=\pi ^{\frak{sp}_{2n}}(\mathtt{s}_{\lambda }^{%
\frak{sp}})\text{ and }s_{\lambda }^{\frak{so}_{N}}=\pi ^{\frak{so}_{N}}(%
\mathtt{s}_{\lambda }^{\frak{s0}}).
\end{equation*}
When $\lambda \notin \mathcal{P}_{N},$ we have $\pi ^{\frak{sp}_{2n}}(%
\mathtt{s}_{\lambda }^{\frak{sp}})=0$ and $\pi ^{\frak{so}_{N}}(\mathtt{s}%
_{\lambda }^{\frak{s0}})=0.\;$The situation is more complicated when $%
\lambda \in \mathcal{P}_{N}$ but $\lambda \notin \mathcal{P}_{n},$ that is
if $d(\lambda )$ the number of parts of $\lambda $ verifies $n<d(\lambda
)\leq N$. In this case one shows by using determinantal identities for the
Weyl characters (see \cite{K} Proposition 2.4.1) that 
\begin{equation*}
\pi ^{\frak{sp}_{2n}}(\mathtt{s}_{\lambda }^{\frak{sp}})=\pm s_{\eta _{C}}^{%
\frak{sp}_{2n}}\text{, }\pi ^{\frak{so}_{2n+1}}(\mathtt{s}_{\lambda }^{\frak{%
so}})=\pm s_{\eta _{B}}^{\frak{so}_{2n+1}}\text{ and }\pi ^{\frak{so}_{2n}}(%
\mathtt{s}_{\lambda }^{\frak{so}})=\pm s_{\eta _{D}}^{\frak{so}_{2n}}
\end{equation*}
where the signs $\pm $ and the partitions $\eta _{B},\eta _{C},\eta _{D}$
are determined by simple combinatorial procedures.\ Note that we have in
this case $\left| \eta _{X}\right| \leq \left| \lambda \right| $ for $X=B,C,D
$.

\bigskip

\noindent We shall also need the following proposition (see \cite{K}
Corollary 2.5.3).

\begin{proposition}
\label{prop_inde}Consider a Lie algebra $\frak{g}$ of type $X_{n}\in
\{B_{n},C_{n},D_{n}\}$. Let $\lambda \in \mathcal{P}_{r}$ and $\mu \in 
\mathcal{P}_{s}$.\ Suppose $n\geq r+s$ and set 
\begin{equation*}
V^{\frak{g}}(\lambda )\otimes V^{\frak{g}}(\mu )=\bigoplus_{\nu \in \mathcal{%
P}_{n}}V^{\frak{g}}(\nu )^{\oplus d_{\lambda ,\mu }^{\nu }}.
\end{equation*}
Then the coefficients $d_{\lambda ,\mu }^{\nu }$ do not depend on the rank $n
$ of $\frak{g}$ neither of its type $B,C$ or $D$.
\end{proposition}

\noindent \textbf{Remark: }The previous proposition follows from the
decompositions 
\begin{equation*}
\mathtt{s}_{\lambda }^{\frak{sp}}\times \mathtt{s}_{\mu }^{\frak{sp}%
}=\sum_{\nu \in \mathcal{P}}d_{\lambda ,\mu }^{\nu }\mathtt{s}_{\nu }^{\frak{%
sp}}\text{ and }\mathtt{s}_{\lambda }^{\frak{so}}\times \mathtt{s}_{\mu }^{%
\frak{so}}=\sum_{\nu \in \mathcal{P}}d_{\lambda ,\mu }^{\nu }\mathtt{s}_{\nu
}^{\frak{so}}
\end{equation*}
for any $\lambda ,\mu \in \mathcal{P}$, in the ring $\Lambda $.

\subsection{Lusztig $q$-analogues}

The $q$-analogue $\mathcal{P}_{q}$ of the Kostant partition function
associated to the root system $R$ of the Lie algebra $\frak{g}$ is defined
by the equality 
\begin{equation*}
\prod_{\alpha \in R^{+}}\dfrac{1}{1-qe^{\alpha }}=\sum_{\beta \in \mathbb{Z}%
^{n}}\mathcal{P}_{q}(\beta )e^{\beta }.
\end{equation*}
Note that $\mathcal{P}_{q}^{B_{n}}(\beta )=0$ if $\beta \notin Q^{+}$.\
Given $\lambda $ and $\mu $ two partitions of $\mathcal{P}_{n},$ the Lusztig 
$q$-analogues of weight multiplicity is the polynomial 
\begin{equation*}
K_{\lambda ,\mu }^{\frak{g}}(q)=\sum_{w\in W^{\frak{g}}}(-1)^{\ell (w)}%
\mathcal{P}_{q}(w\circ \lambda -\mu ).
\end{equation*}
It follows from the Weyl character formula that $K_{\lambda ,\mu }^{\frak{g}%
}(1)$ is equal to the dimension of $V_{\mu }^{\frak{g}}(\lambda ).$

\begin{theorem}
(Lusztig \cite{Lut})

\noindent For any partitions $\lambda ,\mu \in \mathcal{P}_{n},$ the
polynomial $K_{\lambda ,\mu }^{\frak{g}}(q)$ has nonnegative integer
coefficients.
\end{theorem}

\noindent We write 
\begin{equation}
K_{\lambda ,\mu }^{\frak{g}}(q)=\sum_{k\geq 0}K_{\lambda ,\mu }^{\frak{g,}%
k}q^{k}.  \label{def_coef_K}
\end{equation}
Then 
\begin{equation}
K_{\lambda ,\mu }^{\frak{g,}k}(q)=\sum_{w\in W^{\frak{g}}}(-1)^{\ell (w)}%
\mathcal{P}^{k}(w\circ \lambda -\mu )  \label{defLszq}
\end{equation}
where for any $\beta \in \mathbb{Z}^{n}$, $\mathcal{P}^{k}(\beta )$ is the
number of ways of decomposing $\beta $ as a sum of $k$ positive roots.

\noindent \textbf{Remark: }One verifies easily that $K_{\lambda ,\mu }^{%
\frak{g}}(q)\neq 0$ only if $\lambda \geq \mu .\;$Moreover, when $\left| \mu
\right| =\left| \lambda \right| ,$ one has $K_{\lambda ,\mu }^{\frak{g}%
}(q)=K_{\lambda ,\mu }^{\frak{gl}_{n}}(q)$ where $K_{\lambda ,\mu }^{\frak{gl%
}_{n}}(q)$ is the Kostka polynomial associated to $(\lambda ,\mu ),$ i.e.\
the Lusztig $q$-analogue associated to the partitions $\lambda ,\mu $ for
the root system $A_{n-1}.$

\bigskip

\noindent We also introduce the Hall-Littlewood polynomials $Q_{\mu
}^{\prime \frak{g}},$ $\mu \in \mathcal{P}_{n}$ defined by 
\begin{equation*}
Q_{\mu }^{\prime \frak{g}}=\sum_{\lambda \in \mathcal{P}_{n}}K_{\lambda ,\mu
}^{\frak{g}}(q)s_{\lambda }^{\frak{g}}\text{.}
\end{equation*}

\subsection{The symmetric algebra $S(\frak{g})$}

Considered as a $G$-module, $\frak{g}$ is irreducible and we have 
\begin{equation}
\left\{ 
\begin{array}{l}
\frak{so}_{2n+1}\simeq V^{\frak{g}}(1,1)\text{ and }\dim (\frak{so}%
_{2n+1})=n(2n+1) \\ 
\frak{sp}_{2n}\simeq V^{\frak{g}}(2)\text{ and }\dim (\frak{sp}_{2n})=n(2n+1)
\\ 
\frak{so}_{2n}\simeq V^{\frak{g}}(1,1)\text{ and }\dim (\frak{so}%
_{2n})=n(2n-1)
\end{array}
\right. .  \label{girreM}
\end{equation}
Let $S(\frak{g})$ be the symmetric algebra associated to $\frak{g}$ and set 
\begin{equation*}
S(\frak{g})=\bigoplus_{k\geq 0}S^{k}(\frak{g})
\end{equation*}
where $S^{k}(\frak{g})$ is the $k$-th symmetric power of $\frak{g}$. By
Proposition \ref{prop_in_A2n} and (\ref{girreM}), we have 
\begin{equation*}
\frak{g}\simeq V^{\frak{gl}_{N}}(1,1)\downarrow _{\frak{g}}^{\frak{gl}_{N}}%
\text{ for }\frak{g=so}_{N}\text{ and }\frak{g}\simeq V^{\frak{gl}%
_{2n}}(2)\downarrow _{\frak{sp}_{2n}}^{\frak{gl}_{2n}}\text{ for }\frak{g=sp}%
_{2n}.
\end{equation*}
This implies the following isomorphisms 
\begin{equation}
S^{k}(\frak{g})\simeq S^{k}(V^{\frak{gl}_{N}}(1,1))\downarrow _{\frak{g}}^{%
\frak{gl}_{N}}\text{ for }\frak{g=so}_{N}\text{ and }S^{k}(\frak{g})\simeq
S^{k}(V^{\frak{gl}_{2n}}(2))\downarrow _{\frak{sp}_{2n}}^{\frak{gl}_{2n}}%
\text{ for }\frak{g=sp}_{2n}.  \label{SArest}
\end{equation}
for any nonnegative integer $k.$

\begin{example}
By using the Weyl dimension formula (see \cite{GW} page 303), one easily
obtain the decompositions 
\begin{equation*}
S^{2}(V^{\frak{gl}_{N}}(1,1))\simeq V^{\frak{gl}_{N}}(1,1,1,1)\oplus V^{%
\frak{gl}_{N}}(2,2)
\end{equation*}
and 
\begin{equation*}
S^{2}(V^{\frak{gl}_{2n}}(2))\simeq V^{\frak{gl}_{2n}}(4)\oplus V^{\frak{gl}%
_{2n}}(2,2).
\end{equation*}
Hence by \ref{SArest} and Proposition \ref{prop_in_A2n}, this gives 
\begin{equation*}
S^{2}(\frak{g})\simeq V^{\frak{g}}(1,1,1,1)\oplus V^{\frak{g}}(2,2)\oplus V^{%
\frak{g}}(2,0)\oplus V^{\frak{g}}(\emptyset )\text{ for }\frak{g=so}_{N}
\end{equation*}
and 
\begin{equation*}
S^{2}(\frak{sp}_{2n})\simeq V^{\frak{sp}_{2n}}(4)\oplus V^{\frak{sp}%
_{2n}}(2,2)\oplus V^{\frak{sp}_{2n}}(1,1)\oplus V^{\frak{sp}_{2n}}(\emptyset
).
\end{equation*}
\end{example}

\noindent \textbf{Remark: }By the previous formulas, the multiplicities
appearing in the decomposition of the square symmetric power of the Lie
algebra $\frak{g}$ of type $X_{n}\in \{B_{n},C_{n},D_{n}\}$ do not depend on
its rank providing $n\geq 2$. We give in Proposition \ref{prop-dec_charS},
the general explicit decomposition of $S^{k}(\frak{g})$ into its irreducible
components.

\section{Graded characters}

\subsection{Graded character of the symmetric algebra}

\noindent Let $V$ be a $G$ or $GL_{n}$-module.\ For any nonnegative integer $%
k,$ write $S^{k}(V)$ for the $k$-th symmetric power of $V$ and set $%
S(V)=\oplus _{k\geq 0}S^{k}(V).\;$Then $S^{k}(V)$ and $S(V)$ are also $G$%
-modules.\ The graded character of $S(V)$ is defined by 
\begin{equation*}
\mathrm{char}_{q}(S(V))=\sum_{k\geq 0}\mathrm{char}(S^{k}(V))q^{k}.
\end{equation*}
Denote by $\mathcal{W}(V)$ the collection of weight of the module $V$
counted with their multiplicities.\ Then we have 
\begin{equation*}
\mathrm{char}_{q}(S(V))=\prod_{\beta \in \mathcal{W}(V)}\frac{1}{1-qe^{\beta
}}.
\end{equation*}

\bigskip

The weights of the Lie algebra $\frak{g}$ of rank $n$ considered as a $G$%
-module are such that 
\begin{equation*}
\mathcal{W}(\frak{g})=\{\alpha \in R,\underset{n\text{ times}}{\underbrace{%
0,...,0}}\}.
\end{equation*}
Thus the graded character $\mathrm{char}_{q}(S(\frak{g}))$ of $S(\frak{g})$
verifies 
\begin{equation}
\mathrm{char}_{q}(S(\frak{g}))=\frac{1}{(1-q)^{n}}\prod_{\alpha \in R}\frac{1%
}{1-qx^{\alpha }}.  \label{defS(g)}
\end{equation}

\begin{proposition}
\label{prop-dec_charS}For any nonnegative integer $k,$ we have 
\begin{eqnarray*}
\mathrm{char}_{q}(S^{k}(\frak{so}_{N}))=\sum_{\lambda \in \mathcal{P}%
_{n}}\sum_{\nu \in \mathcal{P}_{N}^{(1,1)}(2k)}b_{\nu ,\lambda }^{\frak{so}%
_{N}}s_{\lambda }^{\frak{so}_{N}} \\
\mathrm{char}_{q}(S^{k}(\frak{sp}_{2n}))=\sum_{\lambda \in \mathcal{P}%
_{n}}\sum_{\nu \in \mathcal{P}_{2n}^{(2)}(2k)}b_{\nu ,\lambda }^{\frak{sp}%
_{2n}}s_{\lambda }^{\frak{sp}_{2n}}
\end{eqnarray*}
where $b_{\nu ,\lambda }^{\frak{so}_{N}}$ and $b_{\nu ,\lambda }^{\frak{sp}%
_{2n}}$ are the branching coefficients defined in (\ref{dec_char}).
\end{proposition}

\begin{proof}
Suppose first $\frak{g=sp}_{2n}.$ Recall the classical identity 
\begin{equation*}
\prod_{1\leq i\leq j\leq 2n}\frac{1}{1-x_{i}x_{j}}=\sum_{\nu \in \mathcal{P}%
_{2n}^{(2)}}s_{\nu }^{\frak{gl}_{2n}}
\end{equation*}
due to Littlewood. It immediately implies the decomposition 
\begin{equation*}
\prod_{1\leq i\leq j\leq 2n}\frac{1}{1-qx_{i}x_{j}}=\sum_{\nu \in \mathcal{P}%
_{2n}^{(2)}}q^{\frac{\left| \nu \right| }{2}}s_{\nu }^{\frak{gl}%
_{2n}}=\sum_{k\geq 0}\sum_{\nu \in \mathcal{P}_{2n}^{(2)}(2k)}s_{\nu }^{%
\frak{gl}_{2n}}q^{k}.
\end{equation*}
By applying the restriction map $r^{\frak{sp}_{2n}},$ this gives 
\begin{multline*}
\frac{1}{(1-q)^{n}}\prod_{1\leq i<j\leq n}\frac{1}{1-q\frac{x_{i}}{x_{j}}}%
\frac{1}{1-q\frac{x_{j}}{x_{i}}}\prod_{1\leq r\leq s\leq n}\frac{1}{%
1-qx_{r}x_{s}}\frac{1}{1-q\frac{1}{x_{r}x_{s}}}= \\
\sum_{k\geq 0}\sum_{\nu \in \mathcal{P}_{2n}^{(2)}(2k)}s_{\nu }^{\frak{gl}%
_{2n}}(x_{1},...,x_{n},x_{n}^{-1},...,x_{1}^{-1})q^{k}.
\end{multline*}
From (\ref{dec_char}), this can be rewritten on the form 
\begin{equation*}
\mathrm{char}_{q}(S(\frak{sp}_{2n}))=\frac{1}{(1-q)^{n}}\prod_{\alpha \in R}%
\frac{1}{1-qx^{\alpha }}=\sum_{k\geq 0}\sum_{\lambda \in \mathcal{P}%
_{n}}\sum_{\nu \in \mathcal{P}_{2n}^{(2)}(2k)}b_{\nu ,\lambda }^{\frak{sp}%
_{2n}}s_{\lambda }^{\frak{sp}_{2n}}q^{k}
\end{equation*}
which gives the desired identity by considering the coefficient in $q^{k}$.

\noindent When $\frak{g=so}_{2n+1}$ or $\frak{g=so}_{2n},$ one uses the
identity 
\begin{equation*}
\prod_{1\leq i<j\leq 2n}\frac{1}{1-qx_{i}x_{j}}=\sum_{\nu \in \mathcal{P}%
_{2n}^{(1,1)}}q^{\frac{\left| \nu \right| }{2}}s_{\nu }^{\frak{gl}%
_{2n}}=\sum_{k\geq 0}\sum_{\nu \in \mathcal{P}_{2n}^{(1,1)}(2k)}s_{\nu }^{%
\frak{gl}_{2n}}q^{k}
\end{equation*}
and our result follows by similar arguments.
\end{proof}

\bigskip

\noindent In the sequel, we set 
\begin{equation*}
m_{k,\lambda }^{\frak{so}_{N}}=\sum_{\nu \in \mathcal{P}_{N}^{(1,1)}(2k)}b_{%
\nu ,\lambda }^{\frak{so}_{N}}\text{ and }m_{k,\lambda }^{\frak{sp}%
_{2n}}=\sum_{\nu \in \mathcal{P}_{2n}^{(2)}(2k)}b_{\nu ,\lambda }^{\frak{sp}%
_{2n}}.
\end{equation*}
Thus we have $\mathrm{char}_{q}(S^{k}(\frak{so}_{N}))=\sum_{\lambda \in 
\mathcal{P}_{n}}m_{k,\lambda }^{\frak{so}_{N}}s_{\lambda }^{\frak{so}_{N}}$
and $\mathrm{char}_{q}(S^{k}(\frak{sp}_{2n}))=\sum_{\lambda \in \mathcal{P}%
_{n}}m_{k,\lambda }^{\frak{sp}_{2n}}s_{\lambda }^{\frak{sp}_{2n}}$.

\subsection{Universal graded characters $\mathrm{char}_{q}(\mathtt{S}^{\frak{%
sp}})$ and $\mathrm{char}_{q}(\mathtt{S}^{\frak{so}})$}

\begin{proposition}
\label{prop_stabS}Consider a nonnegative integer $k$ and a partition $%
\lambda \in \mathcal{P}_{m}.$ Suppose $n\geq 2k$.\ Then we have the
identities 
\begin{eqnarray}
m_{k,\lambda }^{\frak{so}_{2n+1}}=m_{k,\lambda }^{\frak{so}_{2n}}=\sum_{\nu
\in \mathcal{P}^{(1,1)}(2k)}\sum_{\gamma \in \mathcal{P}^{(2)}}c_{\lambda
,\gamma }^{\nu }  \label{mstable} \\
m_{k,\lambda }^{\frak{sp}_{2n}}=\sum_{\nu \in \mathcal{P}^{(2)}(2k)}\sum_{%
\gamma \in \mathcal{P}^{(1,1)}}c_{\lambda ,\gamma }^{\nu }.  \notag
\end{eqnarray}
In particular, the multiplicities $m_{k,\lambda }^{\frak{sp}%
_{2n}},m_{k,\lambda }^{\frak{so}_{2n+1}}$ and $m_{k,\lambda }^{\frak{so}%
_{2n}}$ do not depend on $n$.
\end{proposition}

\begin{proof}
For any $\nu \in \mathcal{P}_{N}(2k)$ we have $\nu \in \mathcal{P}_{N}(n)$
since $n\geq 2k.\;$We can thus deduce from Propositions \ref{prop_in_A2n}
and \ref{prop-dec_charS} the decompositions 
\begin{equation*}
m_{k,\lambda }^{\frak{so}_{N}}=\sum_{\nu \in \mathcal{P}_{n}^{(1,1)}(2k)}%
\sum_{\gamma \in \mathcal{P}_{2k}^{(2)}}c_{\lambda ,\gamma }^{\nu }\text{
and }m_{k,\lambda }^{\frak{sp}_{2n}}=\sum_{\nu \in \mathcal{P}%
_{n}^{(2)}(2k)}\sum_{\gamma \in \mathcal{P}_{2k}^{(1,1)}}c_{\lambda ,\gamma
}^{\nu }.
\end{equation*}
Since $c_{\lambda ,\gamma }^{\nu }=0$ when $\left| \lambda \right| +\left|
\gamma \right| \neq 2k,$ $m_{k,\lambda }^{\frak{so}_{N}}$ and $m_{k,\lambda
}^{\frak{sp}_{2n}}$ can be rewritten as in (\ref{mstable}) and thus, do not
depend on $n.$
\end{proof}

\bigskip

\noindent We set 
\begin{equation}
m_{k,\lambda }^{\frak{so}}=\lim_{n\rightarrow \infty }m_{k,\lambda }^{\frak{%
so}_{2n+1}}=\sum_{\nu \in \mathcal{P}^{(1,1)}(2k)}\sum_{\gamma \in \mathcal{P%
}^{(1,1)}}c_{\lambda ,\gamma }^{\nu }\text{ and }m_{k,\lambda }^{\frak{sp}%
}=\lim_{n\rightarrow \infty }m_{k,\lambda }^{\frak{sp}_{2n}}=\sum_{\nu \in 
\mathcal{P}^{(2)}(2k)}\sum_{\gamma \in \mathcal{P}^{(2)}}c_{\lambda ,\gamma
}^{\nu }\text{.}  \label{mul_lim}
\end{equation}

\begin{lemma}
\label{lemma_mCD}For any nonnegative integer $k$ and any partition $\lambda ,
$ we have

\begin{enumerate}
\item  $m_{k,\lambda }^{\frak{sp}}=m_{k,\lambda ^{\prime }}^{\frak{so}}$

\item  $m_{k,\lambda }^{\frak{sp}}=m_{k,\lambda }^{\frak{so}}=0$ if $\left|
\lambda \right| >2k$

\noindent where $m_{k,\lambda }^{\frak{sp}}$ and $m_{k,\lambda ^{\prime }}^{%
\frak{so}}$ are the multiplicities defined in (\ref{mul_lim}).
\end{enumerate}
\end{lemma}

\begin{proof}
Write $\iota $ for the bijective map defined on $\mathcal{P}$ by $\lambda
\longmapsto \lambda ^{\prime }$. We have then $\iota (\mathcal{P}^{(2)})=%
\mathcal{P}^{(1,1)}$ and $\iota (\mathcal{P}^{(2)}(2k))=\mathcal{P}%
^{(1,1)}(2k).$ Moreover $c_{\lambda ,\gamma }^{\nu }=c_{\lambda ^{\prime
},\gamma ^{\prime }}^{\nu ^{\prime }}$ for any partitions $\lambda ,\gamma
,\nu .$ This implies $1$ from the definition (\ref{mul_lim}) of $%
m_{k,\lambda }^{\frak{sp}}$ and $m_{k,\lambda ^{\prime }}^{\frak{so}}.$

\noindent Recall that $c_{\lambda ,\gamma }^{\nu }=0$ when $\left| \nu
\right| \neq \left| \lambda \right| +\left| \gamma \right| .$ Since $\left|
\nu \right| =2k$ in the equalities of (\ref{mul_lim}), we have $m_{k,\lambda
}^{\frak{sp}}=m_{k,\lambda }^{\frak{so}}=0$ if $\left| \lambda \right| >2k.$
\end{proof}

\bigskip

We define the universal graded characters $\mathrm{char}_{q}(\mathtt{S}^{%
\frak{sp}})$ and $\mathrm{char}_{q}(\mathtt{S}^{\frak{so}})$ by setting 
\begin{equation}
\mathrm{char}_{q}(\mathtt{S}(\frak{sp)})=\prod_{1\leq i\leq j}\frac{1}{%
1-qx_{i}x_{j}}\text{ and }\mathrm{char}_{q}(\mathtt{S}(\frak{sp)}%
)=\prod_{1\leq i<j}\frac{1}{1-qx_{i}x_{j}}.  \label{def_gradS}
\end{equation}
Note that $\mathrm{char}_{q}(\mathtt{S}^{\frak{sp}})$ and $\mathrm{char}_{q}(%
\mathtt{S}^{\frak{so}})$ belong to the ring $\mathbb{\Lambda \lbrack \lbrack 
}q\mathbb{]]}$ of formal series with coefficients in $\Lambda .\;$

\noindent For any $F=\sum_{k\geq 0}c_{k}q^{k}$ in $\Lambda \lbrack \lbrack
q]],$ the specialization homomorphisms $\pi ^{\frak{sp}_{2n}},\pi ^{\frak{so}%
_{2n+1}}$ and $\pi ^{\frak{so}_{2n}}$ are then defined by setting 
\begin{equation}
\pi ^{\frak{g}}(F)=\sum_{k\geq 0}\pi ^{\frak{g}}(c_{k})q^{k}.
\label{def_piq}
\end{equation}
By (\ref{defS(g)}) and (\ref{def_gradS}), we have then 
\begin{equation}
\pi ^{\frak{so}_{N}}(\mathrm{char}_{q}(\mathtt{S}(\frak{so)}))=\mathrm{char}%
_{q}(S(\frak{so}_{N}))\text{ and }\pi ^{\frak{sp}_{2n}}(\mathrm{char}_{q}(%
\mathtt{S}(\frak{sp)}))=\mathrm{char}_{q}(S(\frak{sp}_{2n})).
\label{act-pi-S}
\end{equation}
Similarly, the linear involution $\varphi $ (see (\ref{def_fi})) is defined
on $\Lambda \lbrack \lbrack q]]$ by 
\begin{equation}
\varphi (F)=\sum_{k\geq 0}\varphi (c_{k})q^{k}.  \label{def_fiq}
\end{equation}
Observe that $\{q^{k}\mathtt{s}_{\lambda }^{\frak{gl}}\mid k\in \mathbb{N}%
,\lambda \in \mathcal{P}\},$ $\{q^{k}\mathtt{s}_{\lambda }^{\frak{so}}\mid
k\in \mathbb{N},\lambda \in \mathcal{P}\}$ and $\{q^{k}\mathtt{s}_{\lambda
}^{\frak{so}}\mid k\in \mathbb{N},\lambda \in \mathcal{P}\}$ are $\mathbb{Z}$%
-bases of $\Lambda \lbrack \lbrack q]]$.

\begin{proposition}
\label{dec_charUs}We have the decompositions 
\begin{eqnarray}
\mathrm{char}_{q}(\mathtt{S}(\frak{so)})=\sum_{k\geq 0}\sum_{\lambda \in 
\mathcal{P}}m_{k,\lambda }^{\frak{so}}\mathtt{s}_{\lambda }^{\frak{so}}q^{k}%
\text{ and}  \label{formal S} \\
\mathrm{char}_{q}(\mathtt{S}(\frak{sp)})=\sum_{k\geq 0}\sum_{\lambda \in 
\mathcal{P}}m_{k,\lambda }^{\frak{sp}}\mathtt{s}_{\lambda }^{\frak{sp}}q^{k}.
\notag
\end{eqnarray}
\end{proposition}

\begin{proof}
One can write 
\begin{eqnarray*}
\mathrm{char}_{q}(\mathtt{S}(\frak{so)})=\prod_{1\leq i<j}\frac{1}{%
1-qx_{i}x_{j}}=\sum_{k\geq 0}\sum_{\nu \in \mathcal{P}^{(1,1)}[2k]}\mathtt{s}%
_{\nu }^{\frak{gl}}q^{k}=\sum_{k\geq 0}\sum_{\lambda \in \mathcal{P}%
}m_{k,\lambda }^{\frak{so}}\mathtt{s}_{\lambda }^{\frak{so}}q^{k} \\
\mathrm{char}_{q}(\mathtt{S}(\frak{sp)})=\prod_{1\leq i\leq j}\frac{1}{%
1-qx_{i}x_{j}}=\sum_{k\geq 0}\sum_{\nu \in \mathcal{P}^{(2)}[2k]}\mathtt{s}%
_{\nu }^{\frak{gl}}q^{k}=\sum_{k\geq 0}\sum_{\lambda \in \mathcal{P}%
}m_{k,\lambda }^{\frak{sp}}\mathtt{s}_{\lambda }^{\frak{sp}}q^{k}
\end{eqnarray*}
where the rigthmost equalities follow from (\ref{deuSgl}).
\end{proof}

\noindent By using $1$ of Lemma \ref{lemma_mCD}, one derives the following
corollary : 

\begin{corollary}
\label{cor_fiS}We have $\varphi (\mathrm{char}_{q}(\mathtt{S}(\frak{so)}))=%
\mathrm{char}_{q}(\mathtt{S}(\frak{sp)}).$
\end{corollary}

\bigskip

\noindent \textbf{Remark: }By $2$ of Lemma \ref{lemma_mCD}, one has 
\begin{equation}
\mathrm{char}_{q}(\mathtt{S}^{k}(\frak{so)})=\sum_{\lambda \in \mathcal{P}%
[2k]}m_{k,\lambda }^{\frak{so}}\mathtt{s}_{\lambda }^{\frak{so}}\text{ and }%
\mathrm{char}_{q}(\mathtt{S}^{k}(\frak{sp)})=\sum_{\lambda \in \mathcal{P}%
[2k]}m_{k,\lambda }^{\frak{sp}}\mathtt{s}_{\lambda }^{\frak{sp}}
\label{dec_charSk}
\end{equation}
where $\mathcal{P}[2k]$ is the set of partitions $\lambda $ such that $%
\left| \lambda \right| \leq 2k$.

\subsection{Universal graded character for harmonic polynomials}

Let $\frak{g}$ be a Lie algebra of type $X_{n}\in \{B_{n},C_{n},D_{n}\}.\;$%
Since the symmetric algebra $S(\frak{g)}$ can be regarded as a $G$-module,
one can consider 
\begin{equation*}
S(\frak{g})^{G}=\{x\in S(\frak{g)}\mid g\cdot x=x\text{ for any }g\in G\}
\end{equation*}
the ring of the $G$-invariants in $S(\frak{g}).$ By a classical theorem of
Kostant \cite{K} we have 
\begin{equation}
S(\frak{g})=H(\frak{g})\otimes S(\frak{g})^{G}  \label{Kot}
\end{equation}
where $H(\frak{g})$ is the ring of $G$-harmonic polynomials.\ The ring $S(%
\frak{g})^{G}$ is generated by algebraically independent homogeneous
polynomials of degrees $d_{i}=m_{i}+1$ and the graded character of $S(\frak{g%
})^{G}$ considered as a $G$-module verifies 
\begin{equation*}
\mathrm{char}_{q}(S(\frak{g})^{G})=\prod_{i=1}^{n}\frac{1}{1-q^{d_{i}}}.
\end{equation*}
By (\ref{Kot}) the graded character of $H(\frak{g})$ can be written 
\begin{equation}
\mathrm{char}_{q}(H(\frak{g}))=\frac{\mathrm{char}_{q}(S(\frak{g}))}{\mathrm{%
char}_{q}(S(\frak{g})^{G})}=\prod_{i=1}^{n}(1-q^{d_{i}})\ \mathrm{char}%
_{q}(S(\frak{g}))=\sum_{k\geq 0}\mathrm{char}(H^{k}(\frak{g}))q^{k}.
\label{def_Hqfin}
\end{equation}

\bigskip

We define the universal graded characters $\mathrm{char}_{q}(H(\frak{sp}))$
and $\mathrm{char}_{q}(H(\frak{so}))$ by setting 
\begin{eqnarray}
\mathrm{char}_{q}(\mathtt{H}(\frak{so})) &=&\prod_{i\geq 1}(1-q^{2i})\ 
\mathrm{char}_{q}(\mathtt{S}(\frak{so}))=\prod_{i\geq
1}(1-q^{2i})\prod_{1\leq i<j}\frac{1}{1-qx_{i}x_{j}}  \label{def_H} \\
\mathrm{char}_{q}(\mathtt{H}(\frak{sp})) &=&\prod_{i\geq 1}(1-q^{2i})\ 
\mathrm{char}_{q}(\mathtt{S}(\frak{sp}))=\prod_{i\geq
1}(1-q^{2i})\prod_{1\leq i\leq j}\frac{1}{1-qx_{i}x_{j}}  \notag
\end{eqnarray}
The universal characters $\mathrm{char}_{q}(\mathtt{S}(\frak{sp}))$ and $%
\mathrm{char}_{q}(\mathtt{S}(\frak{so}))$ belong to $\Lambda \lbrack \lbrack
q]].\;$Moreover $\prod_{i\geq 1}(1-q^{2i})\in \Lambda \lbrack \lbrack q]]$
since it is a formal series in $q$ with integer coefficients.\ Hence $%
\mathrm{char}_{q}(\mathtt{H}(\frak{sp}))$ and $\mathrm{char}_{q}(\mathtt{H}(%
\frak{so}))$ also belong to $\Lambda \lbrack \lbrack q]].$ We set 
\begin{eqnarray*}
\mathrm{char}_{q}(\mathtt{H}(\frak{sp})) &=&\sum_{k\geq 0}\sum_{\lambda \in 
\mathcal{P}}K_{\lambda ,\emptyset }^{\frak{sp},k}\mathtt{s}_{\lambda }^{%
\frak{sp}}q^{k}\text{ and} \\
\mathrm{char}_{q}(\mathtt{H}(\frak{so})) &=&\sum_{k\geq 0}\sum_{\lambda \in 
\mathcal{P}}K_{\lambda ,\emptyset }^{\frak{so},k}\mathtt{s}_{\lambda }^{%
\frak{so}}q^{k}.
\end{eqnarray*}

\begin{lemma}
For any nonnegative integer $k,$ we have 
\begin{equation}
\mathrm{char}_{q}(\mathtt{H}^{k}(\frak{so)})=\sum_{\lambda \in \mathcal{P}%
[2k]}K_{\lambda ,\emptyset }^{\frak{so},k}\mathtt{s}_{\lambda }^{\frak{so}}%
\text{ and }\mathrm{char}_{q}(\mathtt{H}^{k}(\frak{sp)})=\sum_{\lambda \in 
\mathcal{P}[2k]}K_{\lambda ,\emptyset }^{\frak{sp},k}\mathtt{s}_{\lambda }^{%
\frak{sp}}.  \label{dec_charHk}
\end{equation}
Moreover, for any nonnegative integer $n\geq 2k$ 
\begin{equation}
\pi ^{\frak{so}_{N}}(\mathrm{char}_{q}(\mathtt{H}^{k}(\frak{so)}))=\mathrm{%
char}_{q}(H^{k}(\frak{so}_{N}))\text{ and }\pi ^{\frak{sp}_{2n}}(\mathrm{char%
}_{q}(\mathtt{H}^{k}(\frak{sp)}))=\mathrm{char}_{q}(H^{k}(\frak{sp}_{2n})).
\label{restrict}
\end{equation}
\end{lemma}

\begin{proof}
By definition of the coefficients $K_{\lambda ,\emptyset }^{\frak{so},k},$
we have 
\begin{equation}
\mathrm{char}_{q}(\mathtt{H}^{k}(\frak{so)})=\sum_{\lambda \in \mathcal{P}%
}K_{\lambda ,\emptyset }^{\frak{so},k}\mathtt{s}_{\lambda }^{\frak{so}}.
\label{for2}
\end{equation}
Write $[k/2]$ for the quotient of the Euclidean division of $k$ by $2$.\
Then by (\ref{def_H}), $K_{\lambda ,\emptyset }^{\frak{so},k}$ is the
coefficient of in $q^{k}\mathtt{s}_{\lambda }^{\frak{so}}$ appearing in the
expansion of 
\begin{equation}
\psi _{k}=\prod_{i=1}^{[k/2]}(1-q^{2i})\sum_{a=0}^{k}q^{a}\mathrm{char}_{q}(%
\mathtt{S}^{a}(\frak{so)})  \label{For}
\end{equation}
on the basis $\{q^{k}\mathtt{s}_{\lambda }^{\frak{so}}\mid k\in \mathbb{N}%
,\lambda \in \mathcal{P}\}$.\ Indeed, for any $i>[k/2],$ we have $i>k$.\ By (%
\ref{dec_charSk}), 
\begin{equation*}
\mathrm{char}_{q}(\mathtt{S}^{a}(\frak{so)})=\sum_{\lambda \in \mathcal{P}%
[2a]}m_{a,\lambda }^{\frak{so}}\mathtt{s}_{\lambda }^{\frak{so}}
\end{equation*}
for any nonnegative integer $a.\;$Since the integers $a$ appearing in (\ref
{For}) as lower than $k,$ we have by $2$ of Lemma \ref{lemma_mCD} $%
m_{a,\lambda }^{\frak{so}}=0$ for any partition $\lambda $ such that $\left|
\lambda \right| >2k$.\ Hence, the coefficient of $\psi _{k}$ on $%
q^{b}s_{\lambda }^{\frak{so}}$ is equal to $0$ when $\left| \lambda \right|
>2k.$ This permits to consider only the partitions of $\mathcal{P}[2k]$ in
the sum (\ref{for2}) and yields (\ref{dec_charHk}).\ The proof is the same
for $\mathrm{char}_{q}(\mathtt{H}^{k}(\frak{sp)})$.

\noindent We have seen that $\mathrm{char}_{q}(\mathtt{H}^{k}(\frak{so)}$ is
the coefficient in $q^{k}$ of $\psi _{k}.\;$Now for any rank $n\geq 2k,$ $%
\mathrm{char}_{q}(H^{k}(\frak{so}_{N}))$ is the coefficient of $q^{k}$ in 
\begin{equation*}
\phi _{k}=\prod_{i=1}^{[k/2]}(1-q^{2i})\sum_{a=0}^{k}q^{a}\mathrm{char}%
_{q}(S^{a}(\frak{so}_{N}\frak{)})
\end{equation*}
Indeed we have 
\begin{equation*}
\mathrm{char}_{q}(H(\frak{so}_{N}))=\prod_{i=1}^{n}(1-q^{d_{i}})\ \mathrm{%
char}_{q}(S(\frak{so}_{N}))
\end{equation*}
and for any $i>[k/2],$ $d_{i}>k$. By (\ref{act-pi-S}), we know that $\pi ^{%
\frak{so}_{N}}(\mathrm{char}_{q}(\mathtt{S}^{a}(\frak{so)})=\mathrm{char}%
_{q}(S^{a}(\frak{so}_{N}\frak{)}).$ Thus this gives the equality $\pi ^{%
\frak{so}_{N}}(\psi _{k})=\phi _{k}$.\ Hence we have $\pi ^{\frak{so}_{N}}(%
\mathrm{char}_{q}(\mathtt{H}^{k}(\frak{so)}))=\mathrm{char}_{q}(H^{k}(\frak{%
so}_{N})).$

\noindent The proof is similar for $\mathrm{char}_{q}(\mathtt{H}^{k}(\frak{%
sp)})$.
\end{proof}

\bigskip

\noindent \textbf{Remark: }Since (\ref{restrict}) only holds for $n\geq 2k,$
we have 
\begin{equation*}
\pi ^{\frak{so}_{N}}(\mathrm{char}_{q}(\mathtt{H}(\frak{so)}))\neq \mathrm{%
char}_{q}(H(\frak{so}_{N}))\text{ and }\pi ^{\frak{sp}_{2n}}(\mathrm{char}%
_{q}(\mathtt{H}(\frak{sp)}))\neq \mathrm{char}_{q}(H\frak{sp}_{2n}))
\end{equation*}
i.e. there does not exist identities analogous to (\ref{act-pi-S}) for the $%
\mathrm{char}_{q}(\mathtt{H}(\frak{so)})$ and $\mathrm{char}_{q}(\mathtt{H}(%
\frak{sp)})$.

\section{Stabilization of the coefficients $K_{\protect\lambda ,\protect\mu
}^{\frak{g},k}$}

\subsection{Stabilization of the coefficients $K_{\protect\lambda ,\emptyset
}^{\frak{g},k}$}

The theorem below shows that the coefficients of the expansion of $\mathrm{%
char}_{q}(H(\frak{g}))$ on the basis of Weyl characters are the Lusztig $q$%
-analogues associated to the zero weight (i.e. $\mu =\emptyset $).

\begin{theorem}
\label{Th_Hes}(Hesselink \cite{hes}) We have 
\begin{equation*}
\mathrm{char}_{q}(H(\frak{g}))=\sum_{\lambda \in \mathcal{P}_{n}}K_{\lambda
,\emptyset }^{\frak{g}}(q)s_{\lambda }^{\frak{g}}=\sum_{k\geq
0}\sum_{\lambda \in \mathcal{P}_{n}}K_{\lambda ,\emptyset }^{\frak{g,}%
k}q^{k}s_{\lambda }^{\frak{g}}.
\end{equation*}
The multiplicity of $V^{\frak{g}}(\lambda )$ in the decomposition of $H^{k}(%
\frak{g})$ in its irreducible components is equal to $K_{\lambda ,\emptyset
}^{\frak{g,}k}$.
\end{theorem}

\bigskip

Now fix a nonnegative integer $k$ and choose a rank $n\geq 2k$.\ The
partitions $\lambda $ appearing in (\ref{dec_charHk}) verify $\left| \lambda
\right| \leq 2k$.\ Hence they belong to $\mathcal{P}_{n}$ for $n\geq 2k$.\
Since $\pi ^{\frak{so}_{N}}(s_{\lambda }^{\frak{so}})=s_{\lambda }^{\frak{so}%
_{N}}$ and $\pi ^{\frak{sp}_{2n}}(s_{\lambda }^{\frak{sp}})=s_{\lambda }^{%
\frak{sp}_{2n}}$ for any $\lambda \in \mathcal{P}_{n}$ this gives 
\begin{equation*}
\pi ^{\frak{so}_{N}}(\mathrm{char}_{q}(\mathtt{H}^{k}(\frak{so)}%
))=\sum_{\lambda \in \mathcal{P}_{n}}K_{\lambda ,\emptyset }^{\frak{so}%
,k}s_{\lambda }^{\frak{so}_{N}}\text{ and }\pi ^{\frak{sp}_{2n}}(\mathrm{char%
}_{q}(\mathtt{H}^{k}(\frak{sp)}))=\sum_{\lambda \in \mathcal{P}%
_{n}}K_{\lambda ,\emptyset }^{\frak{sp},k}s_{\lambda }^{\frak{sp}_{2n}}.
\end{equation*}
By using (\ref{restrict}), one obtains
\begin{equation}
\mathrm{char}_{q}(H^{k}(\frak{so}_{N}\frak{)})=\sum_{\lambda \in \mathcal{P}%
_{n}}K_{\lambda ,\emptyset }^{\frak{so},k}s_{\lambda }^{\frak{so}_{N}}\text{
and }\mathrm{char}_{q}(H^{k}(\frak{sp}_{2n}\frak{)})=\sum_{\lambda \in 
\mathcal{P}_{n}}K_{\lambda ,\emptyset }^{\frak{sp},k}s_{\lambda }^{\frak{sp}%
_{2n}}.  \label{dec_ind}
\end{equation}
We can now state the stabilization result for the coefficients $K_{\lambda
,\emptyset }^{\frak{g},k}$.

\begin{proposition}
\label{prop_ine_0}Let $m,k$ be nonnegative integers.\ Consider $\lambda \in 
\mathcal{P}_{m}$ and $\frak{g}$ a Lie algebra of type $X_{n}\in
\{B_{n},C_{n},D_{n}\}$.\ Suppose $n\geq 2k$, then 
\begin{equation*}
K_{\lambda ,\emptyset }^{\frak{sp}_{2n}\frak{,}k}=K_{\lambda ,\emptyset }^{%
\frak{sp,}k}\text{ and }K_{\lambda ,\emptyset }^{\frak{so}_{2n+1}\frak{,}%
k}=K_{\lambda ,\emptyset }^{\frak{so}_{2n},k}=K_{\lambda ,\emptyset }^{\frak{%
so},k}.
\end{equation*}
In particular the coefficients $K_{\lambda ,\emptyset }^{\frak{sp}_{2n}\frak{%
,}k},K_{\lambda ,\emptyset }^{\frak{so}_{2n+1}\frak{,}k}$ and $K_{\lambda
,\emptyset }^{\frak{so}_{2n}\frak{,}k}$ do not depend on the rank $n$.
\end{proposition}

\begin{proof}
By Theorem \ref{Th_Hes} we have 
\begin{equation*}
\mathrm{char}_{q}(H^{k}(\frak{so}_{N}\frak{)})=\sum_{\lambda \in \mathcal{P}%
_{n}}K_{\lambda ,\emptyset }^{\frak{so}_{N},k}s_{\lambda }^{\frak{so}_{N}}%
\text{ and }\mathrm{char}_{q}(H^{k}(\frak{sp}_{2n}\frak{)})=\sum_{\lambda
\in \mathcal{P}_{n}}K_{\lambda ,\emptyset }^{\frak{sp}_{2n},k}s_{\lambda }^{%
\frak{sp}_{2n}}.
\end{equation*}
The Proposition then comes by identifying the coefficients appearing in
these decompositions with those appearing in (\ref{dec_ind}).
\end{proof}

\bigskip

\noindent Thus the coefficients $K_{\lambda ,\emptyset }^{\frak{g,}k}$
depends only on $\lambda ,k$ and the type $X=C$ or $X\in \{B,D\}$ of the Lie
algebra considered when $\mathrm{rank}(\frak{g)}\geq 2k$.\ Proposition \ref
{prop_ine_0} can also be rewritten on the form 
\begin{equation*}
\lim_{n\rightarrow \infty }K_{\lambda ,\emptyset }^{\frak{sp}_{2n}\frak{,}%
k}=K_{\lambda ,\emptyset }^{\frak{sp,}k}\text{ and }\lim_{n\rightarrow
\infty }K_{\lambda ,\emptyset }^{\frak{so}_{N}\frak{,}k}=K_{\lambda
,\emptyset }^{\frak{so,}k}.
\end{equation*}

\subsection{Recurrence formulas for the Lusztig $q$-analogues}

We now recall recurrence formulas established in \cite{lec} and \cite{lec2}
which permit to express the Lusztig $q$-analogues associated to the root
system of type $X_{n}\in \{B_{n},C_{n},D_{n}\}$ in terms of those associated
to the root system of type $X_{n-1}$.\ They can be regarded as
generalizations of the Morris recurrence formula which holds for the Lusztig 
$q$-analogues of type $A$.

\noindent Suppose $\frak{g}$ is of rank $n$ and consider $\nu \in \mathcal{P}%
_{m}$ with $m\leq n.$ For any nonnegative integer $l$, the decomposition of
the $\frak{g}$-module $V^{\frak{g}}(\nu )\otimes V^{\frak{g}}(l)$ into its
irreducible components can be written 
\begin{equation*}
V^{\frak{g}}(\gamma )\otimes V^{\frak{g}}(l)\simeq \bigoplus_{\lambda \in 
\mathcal{P}_{n}}V^{\frak{g}}(\lambda )^{\oplus p_{\gamma ,l}^{\frak{g,}%
\lambda }}.
\end{equation*}
This decomposition can be regarded as the analogue of the Pieri rule in type 
$X_{n}$.

\bigskip

\noindent \textbf{Remark: }The multiplicities $p_{\gamma ,l}^{\frak{g,}%
\lambda }$ are not necessarily equal to $0$ or $1$ as in the original Pieri
rule (i.e. for $\frak{gl}_{n}$). Moreover we can have $p_{\gamma ,l}^{\frak{%
g,}\lambda }\neq 0$ when $\left| \lambda \right| <\left| \gamma \right| +l.$
Nevertheless, when $\left| \lambda \right| =\left| \gamma \right| +l,$ one
shows that the Pieri rule for type $X_{n}$ coincide with the original one.
It means that $p_{\gamma ,l}^{\frak{g,}\lambda }=1$ if $\lambda $ is
obtained by adding an horizontal strip of length $l$ to $\gamma $ (i.e. the $%
l$ boxes of $\lambda /\gamma $ belong to distinct columns) and $p_{\gamma
,l}^{\frak{g,}\lambda }=0$ otherwise.

\noindent The following Lemma is a consequence of Proposition \ref{prop_inde}%
.

\begin{lemma}
\label{lem_util2}Consider $\gamma \in \mathcal{P}_{m}$ and $l\in \mathbb{N}$%
. Suppose that $\frak{g}$ is a Lie algebra of type $X_{n}\in
\{B_{n},C_{n},D_{n}\}$ with $n>m.$ Then the multiplicities $p_{\gamma ,l}^{%
\frak{g,}\lambda }$ are independent of the rank $n$ and the type $X\in
\{B,C,D\}$ of $\frak{g}$.
\end{lemma}

\noindent We set 
\begin{equation*}
p_{\nu ,l}^{\lambda }=\lim_{n\rightarrow \infty }p_{\nu ,l}^{\frak{g,}%
\lambda }.
\end{equation*}

\bigskip

\noindent Let $\mu \in \mathcal{P}_{m}.$ Write $p$ maximal in $\{1,...,m\}$
such that $\nu _{p}-p-\mu _{1}+1\geq 0.$ For any $s\in \{1,...,p\}$ let $%
\gamma (s)$ be the partition of length $m-1$ such that 
\begin{equation}
\gamma (s)=\left\{ 
\begin{tabular}{l}
$(\nu _{2},...,\nu _{m})$ if $s=1$ \\ 
$(\nu _{1}+1,\nu _{2}+1,...,\nu _{s-1}+1,\nu _{s+1},...,\nu _{m})$ if $s\geq
2$%
\end{tabular}
\right. .  \label{def_gmmak}
\end{equation}
Finally set 
\begin{equation}
R_{s}=\nu _{s}-s-\mu _{1}+1.  \label{(Rk)}
\end{equation}

\begin{theorem}
\cite{lec2}\label{Th_mor_expli}

\noindent With the above notation, we have for any partitions $\lambda ,\mu
\in \mathcal{P}_{m}$: 
\begin{align*}
\mathrm{(i)}:K_{\nu ,\mu }^{\frak{so}_{2n+1}}(q)=\sum_{s=1}^{p}(-1)^{s-1}%
\times q^{R_{s}}\times \sum_{r+2a=R_{s}}\sum_{\lambda \in \mathcal{P}%
_{m-1}}p_{\gamma (s),r}^{\frak{so}_{2n-1},\lambda }K_{\lambda ,\mu ^{\flat
}}^{\frak{so}_{2n-1}}(q) \\
\mathrm{(ii)}:K_{\nu ,\mu }^{\frak{sp}_{2n}}(q)=\sum_{s=1}^{p}(-1)^{s-1}%
\times \sum_{r+2a=R_{s}}q^{r+a}\sum_{\lambda \in \mathcal{P}_{m-1}}p_{\gamma
(s),r}^{\frak{sp}_{2n-2},\lambda }K_{\lambda ,\mu ^{\flat }}^{\frak{sp}%
_{2n-2}}(q) \\
\mathrm{(iii)}:K_{\nu ,\mu }^{\frak{so}_{2n}}(q)=\sum_{s=1}^{p}(-1)^{s-1}%
\times q^{R_{s}}\times \sum_{r+2a=R_{s}}\sum_{\lambda \in \mathcal{P}%
_{m-1}}p_{\gamma (s),r}^{\frak{so}_{2n-2},\lambda }K_{\lambda ,\mu ^{\flat
}}^{\frak{so}_{2n-2}}(q)
\end{align*}
where $a\in \mathbb{N}$ and $\mu ^{\flat }=(\mu _{2},...,\mu _{m}).$
\end{theorem}

\noindent By using similar arguments to those given in Example 4 page 243 of 
\cite{mac} , one shows that the polynomials $K_{\lambda ,\mu }^{\frak{g}}(q)$
are monic of degree 
\begin{equation*}
\left\{ 
\begin{tabular}{l}
$\sum_{i=1}^{n}(n-i+1)(\lambda _{i}-\mu _{i})\text{ for }\frak{g=so}_{2n+1}$
\\ 
$\sum_{i=1}^{n}(n-i+1/2)(\lambda _{i}-\mu _{i})\text{ for }\frak{g=sp}_{2n}$
\\ 
$\sum_{i=1}^{n}(n-i)(\lambda _{i}-\mu _{i})\text{ for }\frak{g=so}_{2n}$%
\end{tabular}
\right. .
\end{equation*}
Thank to the recurrence formulas of Theorem \ref{Th_mor_expli}, one can
derive a lower bound for the lowest degree appearing in $K_{\nu ,\mu }^{%
\frak{g}}(q)$ when $K_{\nu ,\mu }^{\frak{g}}(q)\neq 0$.

\begin{proposition}
\label{prop_LB}For any partitions $\nu ,\mu \in \mathcal{P}_{n},$ $K_{\nu }^{%
\frak{g}},_{\mu }(q)=0$ or has degree at least $\frac{\left| \nu \right|
-\left| \mu \right| }{2}.$
\end{proposition}

\begin{proof}
We give the proof for $\frak{g=sp}_{2n},$ the arguments are essentially the
same for $\frak{g=so}_{2n+1}$ and $\frak{g=so}_{2n}.$ We proceed by
induction on $n.\;$For $n=1,$ one has $K_{\nu ,\mu }^{\frak{sp}_{2}}(q)=0$
or $K_{\nu ,\mu }^{\frak{sp}_{2}}(q)=q\frac{^{\left| \lambda \right| -\left|
\mu \right| }}{2}$. Now suppose $K_{\lambda ,\mu ^{\flat }}^{\frak{sp}%
_{2n-2}}(q)$ has degree at least $\frac{\left| \lambda \right| -\left| \mu
^{\flat }\right| }{2}$ or is equal to $0$. Consider $\nu ,\mu \in \mathcal{P}%
_{n}$ such that $K_{\nu ,\mu }^{\frak{sp}_{2n}}(q)\neq 0$.\ We apply
recurrence formula $\mathrm{(ii)}$ of Theorem \ref{Th_mor_expli}. Since $%
K_{\nu ,\mu }^{\frak{sp}_{2n}}(q)\neq 0,$ there exist integers $s$ and $%
r\leq R_{s}$ such that $K_{\lambda ,\mu ^{\flat }}^{\frak{sp}_{2n-2}}(q)\neq
0$ with $p_{\gamma (s),r}^{\frak{sp}_{2n-2},\lambda }\neq 0$. One then have $%
\left| \lambda \right| \geq \left| \gamma (s)\right| -r$. The definition (%
\ref{def_gmmak}) of $\gamma (s)$ gives $\left| \gamma (s)\right| =\left| \nu
\right| -\nu _{s}+s-1.$ By the induction hypothesis, the polynomial $%
q^{r+a}K_{\lambda ,\mu ^{\flat }}^{\frak{sp}_{2n-2}}(q)$ has degree at least 
\begin{equation*}
d=r+a+\frac{\left| \lambda \right| -\left| \mu \right| +\mu _{1}}{2}\geq 
\frac{1}{2}r+a+\frac{\left| \nu \right| -\nu _{s}+s-1-\left| \mu \right|
+\mu _{1}}{2}.
\end{equation*}
One other hand, we have $\frac{1}{2}r+a=\frac{1}{2}R_{s}$ for $r+2a=R_{s}$.\
Recall that $R_{s}=\nu _{s}-s-\mu _{1}+1.$ This finally gives 
\begin{equation*}
d\geq \frac{\nu _{s}-s-\mu _{1}+1}{2}+\frac{\left| \nu \right| -\nu
_{s}+s-1-\left| \mu \right| +\mu _{1}}{2}\geq \frac{\left| \nu \right|
-\left| \mu \right| }{2}.
\end{equation*}
\end{proof}

\bigskip

\noindent \textbf{Remark: }By the previous proposition, the coefficients $%
K_{\nu ,\mu }^{\frak{g},k}$ are all equal to zero when $k<\frac{\left| \nu
\right| -\left| \mu \right| }{2}.$ This implies in particular that we have
the decomposition 
\begin{equation}
Q_{\mu }^{\prime \frak{g}}=\sum_{\nu \in \mathcal{P}_{n}}K_{\nu ,\mu }^{%
\frak{g}}(q)s_{\nu }^{\frak{g}}=\sum_{k\geq 0}\sum_{\nu \in \mathcal{P}%
_{n}[2k+\left| \mu \right| ]}K_{\nu ,\mu }^{\frak{g},k}s_{\nu }^{\frak{g}%
}q^{k}  \label{decQ'}
\end{equation}
for the Hall-Littlewood functions $Q_{\mu }^{\prime \frak{g}}$.

\bigskip

Since the integer $p$ and the partitions $\gamma (s)$ defined above do not
depend on the rank of $\frak{g,}$ we obtain from Lemma \ref{lem_util2} and
Theorem \ref{Th_mor_expli} :

\begin{corollary}
(of Theorem \ref{Th_mor_expli})\label{cor_rec_form}

\noindent For any partitions $\lambda ,\mu \in \mathcal{P}_{m}$ and any
integer $n\geq m$ 
\begin{align*}
\mathrm{(i)}:K_{\nu ,\mu }^{\frak{so}_{2n+1}}(q)=\sum_{s=1}^{p}(-1)^{s-1}%
\times q^{R_{s}}\times \sum_{r+2a=R_{s}}\sum_{\lambda \in \mathcal{P}%
_{m-1}}p_{\gamma (s),r}^{\lambda }K_{\lambda ,\mu ^{\flat }}^{\frak{so}%
_{2n-1}}(q) \\
\mathrm{(ii)}:K_{\nu ,\mu }^{\frak{sp}_{2n}}(q)=\sum_{s=1}^{p}(-1)^{s-1}%
\times \sum_{r+2a=R_{s}}q^{r+a}\sum_{\lambda \in \mathcal{P}_{m-1}}p_{\gamma
(s),r}^{\lambda }K_{\lambda ,\mu ^{\flat }}^{\frak{sp}_{2n-2}}(q) \\
\mathrm{(iii)}:K_{\nu ,\mu }^{\frak{so}_{2n}}(q)=\sum_{s=1}^{p}(-1)^{s-1}%
\times q^{R_{s}}\times \sum_{r+2a=R_{s}}\sum_{\lambda \in \mathcal{P}%
_{m-1}}p_{\gamma (s),r}^{\lambda }K_{\lambda ,\mu ^{\flat }}^{\frak{so}%
_{2n-2}}(q)
\end{align*}
where $a\in \mathbb{N}$ and $\mu ^{\flat }=(\mu _{2},...,\mu _{m}).$
\end{corollary}

\noindent The Lemma below will be useful to derive the recurrence formulas
of paragraph \ref{sub-sec-RF}.

\begin{lemma}
\label{lem_tech4}The partitions $\lambda $ appearing in the right hand sides
of the previous formulas for which there exists a pair $(\gamma (s),r)$ such
that $p_{\gamma (s),r}^{\lambda }\neq 0$ must verify one of the following
assertions:

\begin{enumerate}
\item  $\lambda =\nu $ and then $\mu =\emptyset ,s=1,r=R_{1}=\nu _{1},\gamma
_{1}=\nu ^{\flat },$

\item  $\left| \lambda \right| <\left| \nu \right| ,$

\item  $\left| \lambda \right| =\left| \nu \right| $ with $\lambda \neq \nu $
and then $\mu =\emptyset $ and $\left| \lambda ^{\flat }\right| <\left| \nu
^{\flat }\right| .$
\end{enumerate}

\noindent In particular $\left| \lambda \right| =\left| \nu \right| $ only
if $\mu =\emptyset .$
\end{lemma}

\begin{proof}
Consider $\lambda $ and $(\gamma (s),r)$ such that $p_{\gamma
(s),r}^{\lambda }\neq 0$. We must have $\left| \lambda \right| \leq r+\left|
\gamma (s)\right| =R_{s}-(R_{s}-r)+\left| \gamma (s)\right| .$ By definition
of $R_{s}$ (\ref{(Rk)}) and $\gamma (s)$ (\ref{def_gmmak}) we obtain $\left|
\lambda \right| \leq \left| \nu \right| -\mu _{1}-(R_{s}-r).$ Thus $\left|
\lambda \right| <\left| \nu \right| $ when $\mu \neq \emptyset $ and $\left|
\lambda \right| =\left| \nu \right| $ only if $\mu =\emptyset $ and $r=R_{s}.
$ This permits to restrict ourselves to the case when $\mu =\emptyset ,$ $%
\left| \lambda \right| =\left| \nu \right| $ and $r=R_{s}$.

\noindent Suppose first $\lambda =\nu .\;$Then we must have $s=1$. Otherwise 
$\gamma (s)_{1}=\nu _{1}+1$ and we would have $\lambda _{1}>\nu _{1}.$ This
gives $\mu =\emptyset ,s=1,r=R_{1}=\nu _{1}$ and $\gamma _{1}=\nu ^{\flat }$
as desired.

\noindent Now suppose $\left| \lambda \right| =\left| \nu \right| $ with $%
\lambda \neq \nu $.\ Observe that $\left| \lambda ^{\flat }\right| <\left|
\nu ^{\flat }\right| $ if and only if $\lambda _{1}>\nu _{1}.$ When $s>1,$
we have $\gamma (s)_{1}=\nu _{1}+1>\nu _{1},$ thus $\lambda _{1}>\nu _{1}$
(see Remark before Lemma \ref{lem_util2}). When $s=1,$ we have $\gamma
(1)=(\nu _{2},...,\nu _{n})=\nu ^{\flat }$ and $r=\nu _{1}.$ Since $\left|
\lambda \right| =\left| \nu \right| $ and $p_{\nu ^{\flat },\nu
_{1}}^{\lambda }\neq 0,$ $\lambda $ is obtained by adding an horizontal
strip of length $\nu _{1}$ on $\nu ^{\flat }.\;$This implies $\lambda
_{1}>\nu _{1}$ because the number of columns in $\nu ^{\flat }$ is equal to $%
\nu _{2}\leq \nu _{1}$ and we have assumed $\lambda \neq \nu $.
\end{proof}

\subsection{Stabilization of the coefficients $K_{\protect\lambda ,\protect%
\mu }^{\frak{g},k}$}

\begin{theorem}
\label{Th-Stab}Consider $m$ a nonnegative integer and $\nu ,\mu $ two
partitions such that $\nu \in \mathcal{P}_{m}$ and $\mu \in \mathcal{P}_{a}$%
.\ Let $\frak{g}$ be a Lie algebra of type $X_{n}\in \{B_{n},C_{n},D_{n}\}$
and $k$ a nonnegative integer.$\;$Then for any $n\geq 2k+a$, the
coefficients $K_{\nu ,\mu }^{\frak{g},k}$ do not depend on the rank $n$ of $%
\frak{g}$. Under these hypothesis, we have moreover $K_{\nu ,\mu }^{\frak{so}%
_{2n+1},k}=K_{\nu ,\mu }^{\frak{so}_{2n},k}.$
\end{theorem}

\begin{proof}
Suppose first $\frak{g=so}_{2n+1}.$ We proceed by induction on $a$. Note
that we can suppose $a\leq m,$ otherwise $K_{\nu ,\mu }^{\frak{so}%
_{2n+1},k}=0$ for any rank $n.$ If $a=0,$ then $\mu =\emptyset $ and the
theorem follows directly from Proposition \ref{prop_ine_0}. Suppose now our
theorem true for any partition $\mu ^{\flat }$ of length $a-1$ with $1\leq
a\leq m$ and consider $\mu $ a partition of length $a.$ We then apply the
recurrence formulas of Corollary \ref{cor_rec_form}.\ It follows from (\ref
{(Rk)}) that the integers $R_{s}$ appearing in these formulas do not depend
on the rank $n$ considered. This is also true for the multiplicities $%
p_{\gamma (s),r}^{\lambda }.$

\noindent By our induction hypothesis, for any $p\in \mathbb{N}$, the
coefficients $K_{\lambda ,\mu ^{\flat }}^{\frak{so}_{2n-1},p}$ are
independent on $n$. Indeed $\mu ^{\flat }\in \mathcal{P}_{a-1}$ and $\frak{so%
}_{2n-1}$ has rank $n-1\geq 2k+a-1$.\ The recurrence formulas of Corollary 
\ref{cor_rec_form} imply that each coefficient $K_{\nu ,\mu }^{\frak{so}%
_{2n+1},k}$ can be expressed in terms of the coefficients $K_{\lambda ,\mu
^{\flat }}^{\frak{so}_{2n-1},p}$, the integers $R_{s}$ and $p_{\gamma
(s),r}^{\lambda }$. Moreover, this decomposition is independent of $n$.
Hence $K_{\nu ,\mu }^{\frak{so}_{2n+1},k}$ does not depend on $n.$ By using
similar arguments, we prove that $K_{\nu ,\mu }^{\frak{g},k}$ do not depend
on $n$ when $\frak{g=sp}_{2n}$ or $\frak{so}_{2n}$.

\noindent The equality $K_{\nu ,\mu }^{\frak{so}_{2n+1},k}=K_{\nu ,\mu }^{%
\frak{so}_{2n},k}$ is yet obtained by induction on $a$. It is true for $a=0$
by Proposition \ref{prop_ine_0} and the induction follows from the fact that
the recurrence formulas of Corollary \ref{cor_rec_form} are the same for $%
\frak{so}_{2n+1}$ and $\frak{so}_{2n}$.
\end{proof}

\bigskip

\noindent \textbf{Remark: }The arguments used in the previous proof imply
that it is possible to decompose any Lusztig $q$-analogue $K_{\nu ,\mu }^{%
\frak{g}}(q)$ such that $\mu \neq \emptyset $ in terms of the Lusztig $q$%
-analogues $K_{\lambda ,\emptyset }^{\frak{g}}(q).\;$Moreover this
decomposition is independent on the rank $n$ providing this rank is
sufficiently large.\ In this case the decomposition is the same for $K_{\nu
,\mu }^{\frak{so}_{2n+1}}(q)$ and $K_{\nu ,\mu }^{\frak{so}_{2n}}(q).$
Nevertheless, these two polynomials do not coincide since $K_{\nu ,\emptyset
}^{\frak{so}_{2n+1}}(q)\neq K_{\nu ,\emptyset }^{\frak{so}_{2n}}(q)$ in
general. The previous Theorem also establishes the equality $K_{\nu ,\mu }^{%
\frak{so}_{2n+1,k}}(q)=K_{\nu ,\mu }^{\frak{so}_{2n,k}}(q)$ for any $k\leq 
\frac{n-a}{2}$ where $a$ is the number of nonzero parts in $\mu $.

\bigskip

\noindent By Theorem \ref{Th-Stab}, it makes sense to set 
\begin{equation}
K_{\nu ,\mu }^{\frak{so},k}=\lim_{n\rightarrow \infty }K_{\nu ,\mu }^{\frak{%
so}_{2n+1},k}=\lim_{n\rightarrow \infty }K_{\nu ,\mu }^{\frak{so}_{2n},k}%
\text{ and }K_{\nu ,\mu }^{\frak{sp},k}=\lim_{n\rightarrow \infty }K_{\nu
,\mu }^{\frak{sp}_{2n},k}.  \label{def-Kkgen}
\end{equation}

\subsection{A reformulation in terms of the Brylinski-Kostant filtration}

Recall that the Lusztig $q$-analogue $K_{\lambda ,\mu }^{\frak{g}}(q)$ can
also be characterized from the Brylinski-Kostant filtration on the weight
space $V^{\frak{g}}(\lambda )_{\mu }$ \cite{bry}. Take $e=e_{1}+\cdot \cdot
\cdot +e_{n}\in \frak{u}_{+}$ for a principal nilpotent in $\frak{g}$
compatible with $\frak{h}$.\ The $e$-filtration of $V^{\frak{g}}(\lambda
)_{\mu }$ is the finite filtration $J_{e}(V^{\frak{g}}(\lambda )_{\mu })$
such that 
\begin{equation*}
\{0\}\subset J_{e}^{0}(V^{\frak{g}}(\lambda )_{\mu })\subset J_{e}^{1}(V^{%
\frak{g}}(\lambda )_{\mu })\subset \cdot \cdot \cdot
\end{equation*}
where for any nonnegative integer $k,$ 
\begin{equation*}
J_{e}^{k}(V^{\frak{g}}(\lambda )_{\mu })=\{v\in V^{\frak{g}}(\lambda )_{\mu
}\mid e^{k+1}(v)=0\}.
\end{equation*}
For completeness we also set $J_{e}^{-1}(V^{\frak{g}}(\lambda )_{\mu
})=\{0\}.\;$The following theorem is a consequence of the main result of 
\cite{bry}.

\begin{theorem}
(Brylinski\label{th_bry}) Consider $m$ a nonnegative integer and $\lambda
,\mu \in \mathcal{P}_{m}$.\ Let $\frak{g}$ be a Lie algebra of type $%
X_{n}\in \{B_{n},C_{n}D_{n}\}.\;$Then 
\begin{equation}
K_{\lambda ,\mu }^{\frak{g}}(q)=\sum_{k\geq 0}\dim (J_{e}^{k}(V^{\frak{g}%
}(\lambda )_{\mu })/J_{e}^{k-1}(V^{\frak{g}}(\lambda )_{\mu }))q^{k}.
\label{idenBr}
\end{equation}
\end{theorem}

\noindent By using Theorem \ref{Th-Stab}, the dimension of the space $%
J_{e}^{k}(V^{\frak{g}}(\lambda )_{\mu })$ does not depend on the rank $n$ of 
$\frak{g}$ providing $n$ is sufficiently large. More precisely, we have :

\begin{theorem}
\label{th_stabbry}Consider $\lambda \in \mathcal{P}_{m},$ $\mu \in \mathcal{P%
}_{a}$ and $k\in \mathbb{N}$.\ Let $\frak{g}$ be a Lie algebra of type $%
X_{n}\in \{B_{n},C_{n}D_{n}\}$ with $n\geq 2k+a.$ Then $\dim (J_{e}^{k}(V^{%
\frak{g}}(\lambda )_{\mu }))$ is independent of the rank $n$ of $\frak{g}$.
Moreover we have in this case 
\begin{equation*}
\dim (J_{e}^{k}(V^{\frak{so}_{2n+1}}(\lambda )_{\mu }))=\dim (J_{e}^{k}(V^{%
\frak{so}_{2n}}(\lambda )_{\mu })).
\end{equation*}
\end{theorem}

\begin{proof}
We deduce from Theorem \ref{Th-Stab} and (\ref{idenBr}) that the
coefficients 
\begin{equation*}
K_{\lambda ,\mu }^{X,k}=\dim (J_{e}^{k}(V^{\frak{g}}(\lambda )_{\mu }))-\dim
(J_{e}^{k-1}(V^{\frak{g}}(\lambda )_{\mu }))
\end{equation*}
does not depend on $n$ providing $n\geq 2k+a.$ Under this hypothesis, one
can write
\begin{equation*}
\dim (J_{e}^{k}(V^{\frak{so}_{N}}(\lambda )_{\mu
}))=\sum_{a=0}^{k}K_{\lambda ,\mu }^{\frak{so},a}\text{ and }\dim
(J_{e}^{k}(V^{\frak{sp}_{2n}}(\lambda )_{\mu }))=\sum_{a=0}^{k}K_{\lambda
,\mu }^{\frak{sp},a}
\end{equation*}
Hence $\dim (J_{e}^{k}(V^{\frak{so}_{N}}(\lambda )_{\mu }))$ and $\dim
(J_{e}^{k}(V^{\frak{sp}_{2n}}(\lambda )_{\mu }))$ are independent of $n.$
Moreover, we have 
\begin{equation*}
\dim (J_{e}^{k}(V^{\frak{so}_{2n+1}}(\lambda )_{\mu }))=\dim (J_{e}^{k}(V^{%
\frak{so}_{2n}}(\lambda )_{\mu })).
\end{equation*}
\end{proof}

\bigskip

\noindent \textbf{Remark: }In general, the spaces $J_{e}^{k}(V^{\frak{g}%
}(\lambda )_{\mu }$ depend on the rank $n\geq 2k+a$ considered although
their dimension does not. An interesting problem could consist in the
obtention of explicit bases for the weight spaces $J_{e}^{k}(V^{\frak{g}%
}(\lambda )_{\mu }).$

\section{Limit of Lusztig $q$-analogues}

\subsection{The formal series $K_{\protect\lambda ,\protect\mu }^{\frak{so}%
}(q)$ and $K_{\protect\lambda ,\protect\mu }^{\frak{sp}}(q)$}

The results of Theorem \ref{Th-Stab} suggest to introduce the formal series $%
K_{\lambda ,\mu }^{\frak{so}}(q)$ and $K_{\lambda ,\mu }^{\frak{sp}}(q)$
defined by 
\begin{equation*}
K_{\lambda ,\mu }^{\frak{so}}(q)=\sum_{r\geq 0}K_{\lambda ,\mu }^{\frak{so}%
,k}q^{k}\in \mathbb{N}[[q]]\text{ and }K_{\lambda ,\mu }^{\frak{sp}%
}(q)=\sum_{r\geq 0}K_{\lambda ,\mu }^{\frak{sp},k}q^{k}\in \mathbb{N}[[q]]
\end{equation*}
where the coefficients $K_{\lambda ,\mu }^{\frak{so},k}$ and $K_{\lambda
,\mu }^{\frak{sp},k}$ are those defined in (\ref{def-Kkgen}). Then, $%
K_{\lambda ,\mu }^{\frak{so}}(q)$ and $K_{\lambda ,\mu }^{\frak{sp}}(q)$ can
be regarded as the limits of the Lusztig $q$-analogues $K_{\lambda ,\mu }^{%
\frak{so}_{N}}(q)$ and $K_{\lambda ,\mu }^{\frak{sp}_{2n}}(q)$ when the rang 
$n$ of $\frak{g}$ tends to the infinity.

\noindent We have moreover 
\begin{equation}
\mathrm{char}_{q}(\mathtt{H}(\frak{so}))=\sum_{\lambda \in \mathcal{P}%
}K_{\lambda ,\emptyset }^{\frak{sp}}(q)\mathtt{s}_{\lambda }^{\frak{sp}}%
\text{ and }\mathrm{char}_{q}(\mathtt{H}(\frak{so}))=\sum_{\lambda \in 
\mathcal{P}}K_{\lambda ,\emptyset }^{\frak{so}}(q)\mathtt{s}_{\lambda }^{%
\frak{so}}.  \label{charqHlim}
\end{equation}

\noindent \textbf{Remarks:}

\noindent $\mathrm{(i):}$ When $\frak{g=sp}_{2n}$ or $\frak{so}_{2n},$ $%
K_{\lambda ,\mu }^{\frak{g}}(q)=0$ for any partitions $\lambda ,\mu $ such
that $\left| \lambda \right| -\left| \mu \right| $ is odd.\ Thus for such
partitions we have also $K_{\lambda ,\mu }^{\frak{sp}}(q)=K_{\lambda ,\mu }^{%
\frak{so}}(q)=0.$

\noindent $\mathrm{(ii):}$ Observe that we may have $K_{\lambda ,\mu }^{%
\frak{so}_{2n+1}}(q)\neq 0$ even if $\left| \lambda \right| -\left| \mu
\right| $ is odd. In this case we have thus 
\begin{equation}
\lim_{n\rightarrow \infty }K_{\lambda ,\mu }^{\frak{so}_{2n+1,k}}(q)=0
\label{lim0}
\end{equation}
for any nonnegative integer $k.\;$Take as an example $\lambda =(1)$ and $\mu
=\emptyset $.\ Then $K_{(1),\emptyset }^{\frak{so}_{2n+1}}(q)=q^{n-1}$ for
any rank $n\geq 2$. Thus (\ref{lim0}) is verified for any fixed degree $k.$

\subsection{Recurrence formulas for the series $K_{\protect\lambda ,\protect%
\mu }^{\frak{so}}(q)$ and $K_{\protect\lambda ,\protect\mu }^{\frak{sp}}(q)%
\label{sub-sec-RF}$}

By taking the limit when $n$ tends to the infinity in the formulas of
Corollary \ref{cor_rec_form} (which do not depend on $n$), we obtain the
identities 
\begin{align}
K_{\lambda ,\mu }^{\frak{so}}(q)& =\sum_{s=1}^{p}(-1)^{s-1}\times
q^{R_{s}}\times \sum_{r+2a=R_{s}}\sum_{\lambda \in \mathcal{P}}p_{\gamma
(s),r}^{\lambda }K_{\lambda ,\mu ^{\flat }}^{\frak{so}}(q)
\label{rec_for_limit} \\
K_{\nu ,\mu }^{\frak{sp}}(q)& =\sum_{s=1}^{p}(-1)^{s-1}\times
\sum_{r+2a=R_{s}}q^{r+a}\sum_{\lambda \in }p_{\gamma (s),r}^{\lambda
}K_{\lambda ,\mu ^{\flat }}^{\frak{sp}}(q)  \notag
\end{align}
where $a\in \mathbb{N}$ and $\mu ^{\flat }=(\mu _{2},...,\mu _{m}).$ These
identities yield recurrence formulas for the limit of $q$-analogues.

\noindent To see it, suppose first $\mu \neq \emptyset .$ By Lemma \ref
{lem_tech4}, the formal series $K_{\lambda ,\mu ^{\flat }}^{\frak{so}}(q)$
and $K_{\lambda ,\mu ^{\flat }}^{\frak{sp}}(q)$ appearing in the right hand
sides of (\ref{rec_for_limit}) are such that $\left| \lambda \right| <\left|
\nu \right| .\;$Thus formulas (\ref{rec_for_limit}) permit to express the
series $K_{\lambda ,\mu }^{\frak{so}}(q)$ and $K_{\lambda ,\mu }^{\frak{sp}%
}(q)$ respectively in terms of the series $K_{\lambda ,\mu ^{\flat }}^{\frak{%
so}}(q)$ and $K_{\lambda ,\mu ^{\flat }}^{\frak{sp}}(q)$ with $\left|
\lambda \right| <\left| \nu \right| $.

\noindent Now suppose $\mu =\emptyset .\;$Then by Lemma \ref{lem_tech4}$,$ $%
K_{\lambda ,\emptyset }^{\frak{so}}(q)$ and $K_{\lambda ,\emptyset }^{\frak{%
sp}}(q)$ also appear in the right hand side of (\ref{rec_for_limit}) when $%
\gamma (s)=\gamma (1)=\nu ^{\flat },$ $R_{s}=R_{1}=\nu _{1}$.\ We can write 
\begin{align}
K_{\lambda ,\emptyset }^{\frak{so}}(q)& =\frac{1}{1-q^{\nu _{1}}}\left(
q^{\nu _{1}}\sum_{\underset{a\neq 0}{r+2a=\nu _{1}}}\sum_{\lambda \in 
\mathcal{P}}p_{\nu ^{\flat },r}^{\lambda }K_{\lambda ,\emptyset
}^{X}(q)+\sum_{s=2}^{p-1}(-1)^{s-1}\times
q^{R_{s}}\sum_{r+2a=R_{s}}\sum_{\lambda \in \mathcal{P}}p_{\gamma
(s),r}^{\lambda }K_{\lambda ,\emptyset }^{\frak{so}}(q)\right) 
\label{rec_for_limit_0} \\
& \text{and }  \notag \\
K_{\lambda ,\emptyset }^{\frak{sp}}(q)& =\frac{1}{1-q^{\nu _{1}}}\left(
\sum_{\underset{a\neq 0}{r+2a=\nu _{1}}}q^{r+a}\sum_{\lambda \in \mathcal{P}%
}p_{\nu ^{\flat },r}^{\lambda }K_{\lambda ,\emptyset
}^{C}(q)+\sum_{s=2}^{p}(-1)^{s-1}\times
\sum_{r+2a=R_{s}}q^{r+a}\sum_{\lambda \in \mathcal{P}}p_{\gamma
(s),r}^{\lambda }K_{\lambda ,\emptyset }^{\frak{sp}}(q)\right)   \notag
\end{align}
where the series $K_{\lambda ,\emptyset }^{\frak{so}}(q)$ and $K_{\lambda
,\emptyset }^{\frak{sp}}(q)$ appearing in the right hand sides are such that 
$\left| \lambda \right| <\left| \nu \right| $ or, $\left| \lambda \right|
=\left| \nu \right| $ and $\left| \lambda ^{\flat }\right| <\left| \mu
^{\flat }\right| .\;$Thus formulas (\ref{rec_for_limit_0}) permit to express
the series $K_{\nu ,\emptyset }^{\frak{so}}(q)$ and $K_{\nu ,\emptyset }^{%
\frak{sp}}(q)$ respectively in terms of the series $K_{\lambda ,\emptyset }^{%
\frak{so}}(q)$ and $K_{\lambda ,\emptyset }^{\frak{sp}}(q)$ with $\left|
\lambda \right| <\left| \nu \right| $ or, $\left| \lambda \right| =\left|
\nu \right| $ and $\left| \lambda ^{\flat }\right| <\left| \mu ^{\flat
}\right| .$ Observe that $\left| \lambda \right| +\left| \lambda ^{\flat
}\right| <\left| \nu \right| +\left| \nu ^{\flat }\right| $. Hence one can
compute the series $K_{\nu ,\emptyset }^{\frak{so}}(q)$ and $K_{\nu
,\emptyset }^{\frak{sp}}(q)$ by induction on $\left| \nu \right| +\left| \nu
^{\flat }\right| $ starting from the obvious identity $K_{\emptyset
,\emptyset }^{\frak{so}}(q)=K_{\emptyset ,\emptyset }^{\frak{sp}}(q)=1.$

\bigskip

\noindent It thus follows from the previous arguments that the series $%
K_{\nu ,\mu }^{\frak{so}}(q)$ and $K_{\nu ,\mu }^{\frak{sp}}(q)$ with $\mu
\neq \emptyset $ can be computed by induction on $\left| \nu \right| $ from
the series $K_{\nu ,\emptyset }^{\frak{so}}(q)$ and $K_{\nu ,\emptyset }^{%
\frak{sp}}(q).$ The series $K_{\nu ,\emptyset }^{\frak{so}}(q)$ and $K_{\nu
,\emptyset }^{\frak{sp}}(q)$ being obtained by induction on $\left| \nu
\right| +\left| \nu ^{\flat }\right| $ from $K_{\emptyset ,\emptyset }^{%
\frak{so}}(q)=K_{\emptyset ,\emptyset }^{\frak{sp}}(q)=1.$\ We give in
Proposition \ref{prop_expli} explicit formulas for $K_{\nu ,\emptyset }^{%
\frak{so}}(q)$ and $K_{\nu ,\emptyset }^{\frak{sp}}(q)$ when $\nu $ is a row
or a column partition.

\subsection{A duality between the series $K_{\protect\lambda ,\emptyset }^{%
\frak{so}}(q)$ and $K_{\protect\lambda ^{\prime },\emptyset }^{\frak{sp}%
}(q). $}

\begin{proposition}
\label{th_dual}For any partition $\lambda $ we have the duality 
\begin{equation*}
K_{\lambda ,\emptyset }^{\frak{so}}(q)=K_{\lambda ^{\prime },\emptyset }^{%
\frak{sp}}(q)
\end{equation*}
between the limits of the orthogonal and symplectic Lusztig $q$-analogues
corresponding to the weight $0.$
\end{proposition}

\begin{proof}
We have 
\begin{equation*}
\mathrm{char}_{q}(\mathtt{H}(\frak{so}))=\prod_{i\geq 1}(1-q^{2i})\ \mathrm{%
char}_{q}(\mathtt{S}(\frak{so}))\text{ and }\mathrm{char}_{q}(\mathtt{H}(%
\frak{sp}))=\prod_{i\geq 1}(1-q^{2i})\ \mathrm{char}_{q}(\mathtt{S}(\frak{sp}%
)).
\end{equation*}
Moreover by Corollary \ref{cor_fiS}, $\varphi (\mathrm{char}_{q}(\mathtt{S}%
^{k}(\frak{so})))=\mathrm{char}_{q}(\mathtt{S}^{k}(\frak{sp}))$ for any
nonnegative integer $k.$ This implies the equality 
\begin{equation}
\varphi (\mathrm{char}_{q}(\mathtt{H}^{k}(\frak{so})))=\mathrm{char}_{q}(%
\mathtt{H}^{k}(\frak{sp}))\text{ for any }k\in \mathbb{N}\text{.}
\label{dualH}
\end{equation}
Recall that 
\begin{equation*}
\mathrm{char}_{q}(\mathtt{H}^{k}(\frak{so)})=\sum_{\lambda \in \mathcal{P}%
}K_{\lambda ,\emptyset }^{\frak{so},k}\mathtt{s}_{\lambda }^{\frak{so}}\text{
and }\mathrm{char}_{q}(\mathtt{H}^{k}(\frak{sp)})=\sum_{\lambda \in \mathcal{%
P}}K_{\lambda ,\emptyset }^{\frak{sp},k}\mathtt{s}_{\lambda }^{\frak{sp}}
\end{equation*}
By using (\ref{dualH}), this gives 
\begin{equation*}
\mathrm{char}_{q}(\mathtt{H}^{k}(\frak{sp}))=\sum_{\lambda \in \mathcal{P}%
}K_{\lambda ,\emptyset }^{\frak{sp},k}\mathtt{s}_{\lambda }^{\frak{sp}%
}=\sum_{\lambda \in \mathcal{P}}K_{\lambda ,\emptyset }^{\frak{so},k}\mathtt{%
s}_{\lambda ^{\prime }}^{\frak{sp}}.
\end{equation*}
Since the map 
\begin{equation}
\iota :\left\{ 
\begin{array}{c}
\mathcal{P}\rightarrow \mathcal{P} \\ 
\lambda \longmapsto \lambda ^{\prime }
\end{array}
\right.   \label{bij}
\end{equation}
is bijective, we must have $K_{\lambda ,\emptyset }^{\frak{sp}%
,k}(q)=K_{\lambda ^{\prime },\emptyset }^{\frak{so},k}(q)$ for any
nonnegative integer $k$ which proves the proposition.
\end{proof}

\bigskip

\noindent \textbf{Remark}$\mathrm{:}$ The duality of the previous theorem
does not hold for the Lusztig $q$-analogues, that is $K_{\lambda ,\emptyset
}^{\frak{sp}_{2n}}(q)=K_{\lambda ^{\prime },\emptyset }^{\frak{so}%
_{2n+1}}(q) $ in general. Nevertheless we have $K_{\lambda ,\emptyset }^{%
\frak{sp}_{2n},k}(q)=K_{\lambda ^{\prime },\emptyset }^{\frak{so}_{N},k}(q)$
when $k\leq \frac{n}{2}$ according to Proposition \ref{prop_ine_0}.

\subsection{Some explicit formulas}

\noindent We give below some explicit formulas for the series $K_{\nu ,\mu
}^{\frak{so}}(q)$ and $K_{\nu ,\mu }^{\frak{sp}}(q)$ when $\nu $ is a column
or a row partition.\ Note that we have not found such simple formulas for
the Lusztig $q$-analogues $K_{\nu ,\mu }^{\frak{g}}(q)$ even in the case
when $\nu $ is a row or a column.

\begin{proposition}
\label{prop_expli}Consider $l$ a nonnegative integer. Recall that $(2l)$ and 
$(1^{2l})$ are the row and column partitions of length and height $2l,$
respectively. We have 
\begin{eqnarray*}
K_{(2l),\emptyset }^{\frak{sp}}(q) &=&K_{(1^{2l}),\emptyset }^{\frak{so}}(q)=%
\frac{q^{l}}{\prod_{i=1}^{l}(1-q^{2i})} \\
K_{(2l),\emptyset }^{\frak{so}}(q) &=&K_{(1^{2l}),\emptyset }^{\frak{sp}}(q)=%
\frac{q^{2l}}{\prod_{i=1}^{l}(1-q^{2i})}
\end{eqnarray*}
\end{proposition}

\begin{proof}
We only give the proof for the first equality of the proposition. The proof
for the second is similar.

\noindent We use the recurrence formula (\ref{rec_for_limit_0}). We have
then $p=1,R_{1}=2l$ and $\gamma (1)=\emptyset $.\ Thus $p_{\gamma
(1),r}^{\lambda }\neq 0$ only when $\lambda =r$ and in this case $p_{\gamma
(1),r}^{r}=1$.\ This yields for any $l\geq 1$%
\begin{equation}
K_{(2l),\emptyset }^{\frak{sp}}(q)=\frac{1}{1-q^{2l}}\sum_{\underset{r\neq 2l%
}{r+2a=2l}}q^{r+a}K_{(r),\emptyset }^{\frak{sp}}(q)=\frac{q^{l}}{1-q^{2l}}%
\sum_{b=0}^{l-1}q^{b}K_{(2b),\emptyset }^{\frak{sp}}(q)  \label{rec0}
\end{equation}
where the last equality is obtained by setting $r=2b.\;$By an immediate
induction starting from $K_{\emptyset ,\emptyset }^{\frak{sp}}(q)=1,$ one
derives the desired formula 
\begin{equation*}
K_{(2l),\emptyset }^{\frak{sp}}(q)=\frac{q^{l}}{\prod_{i=1}^{l}(1-q^{2i})}
\end{equation*}
by using the identity 
\begin{equation}
\sum_{b=0}^{l-1}q^{b}K_{(2b),\emptyset }^{\frak{sp}}(q)=\sum_{b=0}^{l-1}%
\frac{q^{2b}}{\prod_{i=1}^{b}(1-q^{2i})}=\frac{1}{\prod_{i=1}^{l-1}(1-q^{2i})%
}.  \label{iden}
\end{equation}
We deduce then $K_{(2l),\emptyset }^{\frak{sp}}(q)=K_{(1^{2l}),\emptyset }^{%
\frak{so}}(q)$ from Theorem \ref{th_dual}.
\end{proof}

\begin{corollary}
Consider $m$ a nonnegative integer and $\mu $ a partition with $d$ nonzero
parts. Then

\begin{enumerate}
\item  $K_{(m),\mu }^{\frak{sp}}(q)\neq 0$ and $K_{(m),\mu }^{\frak{so}%
}(q)\neq 0$ only if $m-\left| \mu \right| \in 2\mathbb{N}$.\ In this case 
\begin{equation*}
K_{(m),\mu }^{\frak{sp}}(q)=q^{h(\mu )}K_{(2l),\emptyset }^{\frak{sp}}(q)=%
\frac{q^{h(\mu )+l}}{\prod_{i=1}^{l}(1-q^{2i})}\text{ and }K_{(m),\mu }^{%
\frak{so}}(q)=q^{h(\mu )}K_{(2l),\emptyset }^{\frak{so}}(q)=\frac{q^{h(\mu
)+2l}}{\prod_{i=1}^{l}(1-q^{2i})}
\end{equation*}
where $h(\mu )=\sum_{1\leq i\leq d}(i-1)\mu _{i}$ and $l=\frac{m-\left| \mu
\right| }{2}.$

\item  $K_{(1^{m}),\mu }^{\frak{so}}(q)\neq 0$ and $K_{(m),\mu }^{\frak{so}%
}(q)\neq 0$ only if $\mu =(1^{p})$ with $m-p\in 2\mathbb{N}$ and in this
case 
\begin{equation*}
K_{(1^{m}),(1^{p})}^{\frak{sp}}(q)=K_{(1^{2l}),\emptyset }^{\frak{sp}}(q)=%
\frac{q^{2l}}{\prod_{i=1}^{l}(1-q^{2i})}\text{ and }K_{(1^{m}),(1^{p})}^{%
\frak{so}}(q)=K_{(1^{2l}),\emptyset }^{\frak{so}}(q)=\frac{q^{l}}{%
\prod_{i=1}^{l}(1-q^{2i})}
\end{equation*}
where $l=\frac{m-\left| p\right| }{2}.$
\end{enumerate}
\end{corollary}

\begin{proof}
$1:$ We proceed by induction on $d$ the number of nonzero parts of $\mu $.\
If $d=0,$ the result follows from Proposition \ref{prop_expli}. Suppose $b>0$
and apply the recurrence formula (\ref{rec_for_limit}). We have $p=1,$ $%
R_{1}=m-\mu _{1}$ and $\gamma (1)=\emptyset $.\ This gives 
\begin{equation*}
K_{(m),\mu }^{\frak{sp}}(q)=\sum_{r+2a=m-\mu _{1}}q^{r+a}K_{(r),\mu ^{\flat
}}^{\frak{sp}}(q).
\end{equation*}
Since $K_{(r),\mu ^{\flat }}^{\frak{sp}}(q)=0$ when $r<\left| \mu ^{\flat
}\right| ,$ we can suppose $r\geq \left| \mu ^{\flat }\right| $ in the
previous sum.\ This gives by using the induction hypothesis 
\begin{equation*}
K_{(m),\mu }^{\frak{sp}}(q)=q^{h(\mu ^{\flat })}\sum_{r+2a=m-\mu
_{1}}q^{r+a}K_{r-\left| \mu ^{\flat }\right| ,\emptyset }^{\frak{sp}}(q).
\end{equation*}
We must have $r-\left| \mu ^{\flat }\right| \in 2\mathbb{N}$, thus we can
set $b=\frac{r-\left| \mu ^{\flat }\right| }{2}$. One then obtains 
\begin{equation*}
K_{(m),\mu }^{\frak{sp}}(q)=q^{h(\mu ^{\flat })}\sum_{b=0}^{l}q^{\frac{%
\left| \mu \right| +m}{2}+b-\mu _{1}}K_{2b,\emptyset }^{\frak{sp}%
}(q)=q^{h(\mu )}\sum_{b=0}^{l}q^{l+b}K_{2b,\emptyset }^{\frak{sp}}(q)
\end{equation*}
where the last equality follows from the equalities $l=\frac{m-\left| \mu
\right| }{2}$ and $h(\mu )=h(\mu ^{\flat })+\left| \mu \right| -\mu _{1}$.
By using (\ref{iden}), this gives 
\begin{equation*}
K_{(m),\mu }^{\frak{sp}}(q)=q^{h(\mu )}\times \frac{q^{l}}{%
\prod_{i=1}^{l}(1-q^{2i})}=q^{h(\mu )}K_{(2l),\emptyset }^{\frak{sp}}(q).
\end{equation*}
The proof is similar for $K_{(m),\mu }^{\frak{so}}(q)$.

\noindent $2:$ By applying (\ref{rec_for_limit}), we obtain this time $p=1,$ 
$R_{1}=0$ and $\gamma (1)=\emptyset $. Hence 
\begin{equation*}
K_{(1^{m}),(1^{p})}^{\frak{sp}}(q)=K_{(1^{m-1}),(1^{p-1})}^{\frak{sp}}(q)%
\text{.}
\end{equation*}
By an immediate induction, this gives $K_{(1^{m}),(1^{p})}^{\frak{sp}%
}(q)=K_{(1^{m-p}),\emptyset }^{\frak{sp}}(q)=K_{(1^{2l}),\emptyset }^{\frak{%
sp}}(q)$ and our formula follows from Proposition \ref{prop_expli}.
\end{proof}

\section{Hall-Littlewood polynomials in infinitely many variables}

\subsection{The ring of graded universal characters}

We now consider the ring $\Delta $ generated over $\mathbb{Z[[}q]]$ by the
formal characters $\mathtt{s}_{\lambda }^{\frak{gl}},$ $\lambda \in \mathcal{%
P}$ with multiplication defined by 
\begin{equation*}
\mathtt{s}_{\lambda }^{\frak{gl}}\cdot \mathtt{s}_{\mu }^{\frak{gl}%
}=\sum_{\nu \in \mathcal{P}}c_{\lambda ,\mu }^{\nu }\mathtt{s}_{\nu }^{\frak{%
gl}}
\end{equation*}
and for any pair $F=\sum_{\lambda \in \mathcal{P}}C_{\lambda }\mathtt{s}%
_{\lambda }^{\frak{gl}}$, $G=\sum_{\mu \in \mathcal{P}}C_{\mu }\mathtt{s}%
_{\mu }^{\frak{gl}}$%
\begin{equation*}
F\cdot G=\sum_{\nu \in \mathcal{P}}\sum_{\lambda ,\mu \in \mathcal{P}%
}C_{\lambda }C_{\mu }c_{\lambda ,\mu }^{\nu }\mathtt{s}_{\nu }^{\frak{gl}}.
\end{equation*}
Observe that $F\cdot G$ is well defined since we have $c_{\lambda ,\mu
}^{\nu }=0$ if $\left| \nu \right| \neq \left| \lambda \right| +\left| \mu
\right| $ and thus $\sum_{\lambda ,\mu \in \mathcal{P}}C_{\lambda }C_{\mu
}c_{\lambda ,\mu }^{\nu }$ is finite. Then $\mathcal{B}^{\frak{gl}}=\{%
\mathtt{s}_{\lambda }^{\frak{gl}}\mid \lambda \in \mathcal{P}\},$ is a $%
\mathbb{Z[}[q]]$-basis of $\Delta $.\ We then defined the $\mathbb{Z}[[q]]$%
-bases $\mathcal{B}^{\frak{so}}=\{\mathtt{s}_{\lambda }^{\frak{so}}\mid
\lambda \in \mathcal{P}\}$ and $\mathcal{B}^{\frak{sp}}=\{\mathtt{s}%
_{\lambda }^{\frak{sp}}\mid \lambda \in \mathcal{P}\}$ so that (\ref{deuSgl}%
) is verified. We then write $<\cdot ,\cdot >^{\frak{so}}$ and $<\cdot
,\cdot >^{\frak{sp}}$ respectively for the inner scalar products which make
the bases $\mathcal{B}^{\frak{so}}$ and $\mathcal{B}^{\frak{sp}}$
orthonormal.

\subsection{Hall-Littlewood polynomials in infinitely many variables}

For any partition $\mu $, we set 
\begin{equation*}
\mathtt{Q}_{\mu }^{\prime \frak{so}}=\sum_{\lambda \in \mathcal{P}%
}K_{\lambda ,\mu }^{\frak{so}}(q)\mathtt{s}_{\lambda }^{\frak{so}}\text{ and 
}\mathtt{Q}_{\mu }^{\prime \frak{sp}}=\sum_{\lambda \in \mathcal{P}%
}K_{\lambda ,\mu }^{\frak{sp}}(q)\mathtt{s}_{\lambda }^{\frak{sp}}.
\end{equation*}
Since $K_{\lambda ,\mu }^{\frak{so}}(q)=K_{\lambda ,\mu }^{\frak{sp}}(q)=0$
when $\lambda <\mu $ for the usual order on partitions, the transition
matrices of the families $\{\mathtt{Q}_{\lambda }^{\prime \frak{so}}\mid
\lambda \in \mathcal{P}\}$ and $\{\mathtt{Q}_{\lambda }^{\prime \frak{sp}%
}\mid \lambda \in \mathcal{P}\}$ respectively on the bases $\mathcal{B}^{%
\frak{so}}$ and $\mathcal{B}^{\frak{sp}}$ are upper-unitriangular.\ Hence $\{%
\mathtt{Q}_{\lambda }^{\prime \frak{so}}\mid \lambda \in \mathcal{P}\}$ and $%
\{\mathtt{Q}_{\lambda }^{\prime \frak{sp}}\mid \lambda \in \mathcal{P}\}$
are bases of $\Delta $.\ 

\bigskip

\noindent We then define the basis $\{\mathtt{P}_{\lambda }^{\frak{so}}\mid
\lambda \in \mathcal{P}\}$ (resp.\ $\{\mathtt{P}_{\lambda }^{\frak{sp}}\mid
\lambda \in \mathcal{P})$ as the dual basis of $\{\mathtt{Q}_{\lambda
}^{\prime \frak{so}}\mid \lambda \in \mathcal{P}\}$ (resp.\ $\{\mathtt{Q}%
_{\lambda }^{\prime \frak{sp}}\mid \lambda \in \mathcal{P}\})$ with respect
to $<\cdot ,\cdot >^{\frak{so}}$ (resp.\ $<\cdot ,\cdot >^{\frak{sp}})$.\
Then $\mathtt{P}_{\lambda }^{\frak{so}}$ and $\mathtt{P}_{\lambda }^{\frak{sp%
}}$ can be regarded as Hall-Littlewood polynomials in infinitely many
variables. One has the identities 
\begin{equation*}
\mathtt{s}_{\lambda }^{\frak{so}}=\sum_{\lambda \in \mathcal{P}}K_{\lambda
,\mu }^{\frak{so}}(q)\mathtt{P}_{\mu }^{\frak{so}}\text{ and }\mathtt{s}%
_{\lambda }^{\frak{sp}}=\sum_{\lambda \in \mathcal{P}}K_{\lambda ,\mu }^{%
\frak{sp}}(q)\mathtt{P}_{\mu }^{\frak{sp}}.
\end{equation*}

\noindent \textbf{Remark}$\mathrm{:}$

\noindent $\mathrm{(i):}$\textrm{\ }One cannot obtain Hall-Littelwood
polynomials $\mathtt{Q}_{\lambda }^{\prime \frak{g}}$ in infinitely many
variables by considering the limit when $n$ tends to the infinity in (\ref
{def_coef_K}) and (\ref{defLszq}). Indeed, the number of ways of decomposing
a weight $\beta $ as a sum of $k$ positive roots (where $k$ is a fixed
nonnegative integer) may strictly increase with $n$.\ As an example for $k=2$
we have $2\varepsilon _{1}=(\varepsilon _{1}-\varepsilon _{i})+(\varepsilon
_{1}+\varepsilon _{i})$ for any $i\in \{2,...,n\}$.

\noindent $\mathrm{(ii):}$\textrm{\ }Recall that  the Hall-Littlewood
polynomial $P_{\mu }^{\frak{g}}$ (see \cite{NR}) associated to the partition 
$\mu $ is defined by 
\begin{equation}
P_{\mu }^{\frak{g}}=\frac{1}{W_{\mu }^{\frak{g}}(q)}\sum_{w\in W^{\frak{g}%
}}w\left( e^{\mu }\prod_{\alpha \in R_{+}}\frac{1-qe^{-\alpha }}{%
1-e^{-\alpha }}\right)   \label{def_Pmu}
\end{equation}
where $W_{\mu }^{\frak{g}}(q)=\sum_{w\in W_{\mu }^{\frak{g}}}q^{\ell (w)}$
with $W_{\mu }^{\frak{g}}$ the stabilizer of $\mu $ in $W^{\frak{g}}.$ Since
the number of elements of $W^{\frak{g}}$ of fixed length strictly increase
with the rank of $\frak{g},$ the polynomials $W_{\mu }^{\frak{g}}(q)$ have
no limit in $\mathbb{Z[}[q]]$. This implies that is not possible to define
Hall-Littlewood polynomials in infinitely many variables by taking the limit
when $n$ tends to the infinity in (\ref{def_Pmu}).

\noindent $\mathrm{(iii):}$\textrm{\ }Similarly to (\ref{decQ'}), we have
the decompositions
\begin{eqnarray*}
\mathtt{Q}_{\lambda }^{\prime \frak{so}} &=&\sum_{\nu \in \mathcal{P}}K_{\nu
,\mu }^{\frak{so}}(q)\mathtt{s}_{\nu }^{\frak{so}}=\sum_{k\geq 0}\sum_{\nu
\in \mathcal{P}[2k+\left| \mu \right| ]}K_{\nu ,\mu }^{\frak{so},k}\mathtt{s}%
_{\nu }^{\frak{so}}q^{k}\text{ and} \\
\mathtt{Q}_{\lambda }^{\prime \frak{sp}} &=&\sum_{\nu \in \mathcal{P}}K_{\nu
,\mu }^{\frak{sp}}(q)\mathtt{s}_{\nu }^{\frak{sp}}=\sum_{k\geq 0}\sum_{\nu
\in \mathcal{P}[2k+\left| \mu \right| ]}K_{\nu ,\mu }^{\frak{sp},k}\mathtt{s}%
_{\nu }^{\frak{sp}}q^{k}.
\end{eqnarray*}


\begin{thebibliography}{99}
\bibitem{bry}  \textsc{R-K. Brylinski,} \textit{Limits of weight spaces,
Lusztig's }$q$\textit{-analogs and fiberings of adjoint orbits}, J.\ Amer.\
Math.\ Soc, \textbf{2}, no.3 (1989), 517-533.

\bibitem{FH}  \textsc{W.\ Fulton, J. Harris, }\textit{Representation theory}%
, Graduate Texts in Mathematics, Springer-Verlag.

\bibitem{GW}  \textsc{G. Goodman, N. R Wallach, }\textit{Representation
theory and invariants of the classical groups}, Cambridge University Press.

\bibitem{hes}  \textsc{W-H. Hesselink,} \textit{Characters of the nullcone, }%
Math.\ Ann.\ \textbf{252}, 179-182 (1980).

\bibitem{KN}  \textsc{M. Kashiwara, T. Nakashima,} \textit{Crystal graphs
for representations of the }$q$\textit{-analogue of classical Lie algebras},
Journal of Algebra, \textbf{165}, 295-345 (1994).

\bibitem{K}  \textsc{K. Koike,}\textit{\ }\textsc{I. Terada,}\textit{\ Young
diagrammatic methods for the representations theory of the classical groups
of type }$B_{n},C_{n}$ and $D_{n},$ Journal of Algebra, \textbf{107},
466-511 (1987).

\bibitem{KT}  \textsc{K. Koike,}\textit{\ }\textsc{I. Terada,}\textit{\
Young diagrammatic methods for the restriction of representations of complex
classical Lie groups to reductive subgroups of maximal rank, Advances in
Mathematics, 79, 104-135 (1990).}

\bibitem{kost}  \textsc{B.\ Kostant}, \textit{Lie groups representations on
polynomial rings}, Amer.\ J.\ Math, \textbf{85}, 327-404 (1963).

\bibitem{Li}  \textsc{D-E. Littlewood,} \textit{The theory of group
characters and matrix representations of groups, }Oxford University Press,
second edition (1958).

\bibitem{LS}  \textsc{A. Lascoux, M-P. Sch\"{u}tzenberger, }\ \textit{Le mono%
}$\mathit{\ddot{\imath}}$\textit{de plaxique}, in non commutative structures
in algebra and geometric combinatorics A. de Luca Ed., Quaderni della
Ricerca Scientifica del C.N.R., Roma, 1981\textit{.}

\bibitem{LSc1}  \textsc{A. Lascoux, M-P. Sch\"{u}tzenberger, }\ \textit{Sur
une conjecture de H.O Foulkes}, CR Acad Sci Paris, \textbf{288}, 95-98
(1979).

\bibitem{lec}  \textsc{C.\ Lecouvey,}\textit{\ Kostka-Foulkes polynomials
cyclage graphs and charge statistic for the root system }$C_{n},$ Journal of
Algebraic Combinatorics, \textbf{21}, 203-240 (2005).

\bibitem{lec2}  \textsc{C.\ Lecouvey,}\textit{\ Combinatorics of crystal
graphs and Kostka-Foulkes polynomials for the root systems }$B_{n},C_{n}$%
\textit{\ and }$D_{n},$ European Journal of Combinatorics, \textbf{27},
526-557 (2006).

\bibitem{Lut}  \textsc{G. Lusztig,} \textit{Singularities, character
formulas, and a }$q$\textit{-analog of weight multiplicities}, Analyse et
topologie sur les espaces singuliers (II-III), Asterisque \textbf{101-102},
208-227 (1983).

\bibitem{mac}  \textsc{I-G. Macdonald,} \textit{Symmetric functions and Hall
polynomials}, Second edition, Oxford Mathematical Monograph, Oxford
University Press, New York, (1995).

\bibitem{Mor}  \textsc{A-O. Morris,} \textit{The characters of the group }$%
GL(n,q),$ Math. Zeitschr. \textbf{81}, 112-123 (1963).

\bibitem{NR}  \textsc{K. Nelsen, A. Ram,} \textit{Kostka-Foulkes polynomials
and Macdonald spherical functions}, Surveys In Combinatorics 2003, C.
Wensley ed.\ London Math.\ Soc.\ Lect.\ Notes \textbf{307}, Cambridge
University Press, 325-370, (2003).
\end{thebibliography}
\end{document}